\documentclass[11pt,reqno]{amsart}
\usepackage{amscd, amsthm, amssymb, xypic, multirow, mathtools}

\newtheorem{theorem}{Theorem}[section]
\newtheorem{lemma}[theorem]{Lemma}
\newtheorem{prop}[theorem]{Proposition}

\theoremstyle{definition}
\newtheorem{defi}[theorem]{Definition}
\newtheorem{example}[theorem]{Example}
\theoremstyle{remark}
\newtheorem*{remark}{Remark}
\newtheorem*{remarks}{Remarks}

\numberwithin{equation}{section}

\DeclareMathOperator{\Hom}{Hom}

\DeclareMathOperator{\im}{im}

\DeclareFontFamily{OT1}{rsfs}{}
\DeclareFontShape{OT1}{rsfs}{n}{it}{<-> rsfs10}{}
\DeclareMathAlphabet{\curly}{OT1}{rsfs}{n}{it}
\newcommand{\KK}{\curly{K}\!\text{\it 3}\,}

\newcommand{\mbz}{\mathbb Z}
\newcommand{\mbc}{\mathbb C}
\newcommand{\mbr}{\mathbb{R}}
\newcommand{\mco}{\mathcal{O}}
\newcommand{\CP}{\mathbb{C}P}
\newcommand{\CX}{\mathbb{C}}
\newcommand{\RE}{\mathbb{R}}
\newcommand{\Z}{\mathbb{Z}}
\newcommand{\Q}{\mathbb{Q}}

\newcommand{\ra}{\rightarrow}
\newcommand{\rk}{\operatorname{rk}}
\newcommand{\bfx}{\mathbf{x}}
\newcommand{\oW}{{\overline{W}}}
\newcommand{\we}{\wedge}
\newcommand{\tka}{\tilde{\kappa}}
\renewcommand{\phi}{\varphi}
\renewcommand{\setminus}{-}

\begin{document}

\title[$K3$ surfaces and $G_2$-manifolds]{$K3$ surfaces with non-symplectic
  involution and compact irreducible $G_2$-manifolds}

\author[A. Kovalev]{Alexei Kovalev}
\address{DPMMS, University of Cambridge, Centre for Mathematical
  Sciences,\break
  Wilberforce Road, Cambridge CB3 0WB, England}
\email{a.kovalev@dpmms.cam.ac.uk}

\author[N.-H.\ Lee]{Nam-Hoon Lee$^*$}
\thanks{$^*$Supported by National Research Foundation of Korea
Grant funded by the Korean Government (2010-0015242) and 2010 Hongik
University Research Fund.} 
\address{Department of Mathematics Education, Hongik University 42-1,
Sangsu-Dong, Mapo-Gu, Seoul 121-791, Korea}
\email{nhlee@hongik.ac.kr}

\begin{abstract}
We consider the connected-sum method of constructing compact Riemannian
7-manifolds with holonomy $G_2$ developed by the first named author. The
method requires pairs of projective complex threefolds endowed with
anticanonical $K3$ divisors and the latter $K3$ surfaces should satisfy a
certain `matching condition' intertwining on their periods and K\"ahler
classes. Suitable examples of threefolds were previously obtained
by blowing up curves in Fano threefolds.

In this paper, we give a large new class of suitable algebraic threefolds
using theory of $K3$ surfaces with non-symplectic involution due to Nikulin.
These threefolds are {\em not} obtainable from Fano threefolds as above,
and admit matching pairs leading to topologically new examples of compact
irreducible $G_2$-manifolds. `Geography' of the values of Betti numbers
$b^2$,$b^3$ for the new (and previously known) examples of irreducible
$G_2$ manifolds is also discussed.
\end{abstract}

\maketitle

\section{Introduction}
The $G_2$-manifolds may be defined as 7-dimensional real manifolds endowed
with a positive 3-form which induces a $G_2$-structure, and hence a Riemannian
metric, and is parallel with respect to this metric (see \S\ref{sum} for
details and further references). The holonomy of the metric on
$G_2$-manifold reduces to a subgroup of the exceptional Lie group~$G_2$.
Joyce~\cite{Jo1,Jo2} constructed the first examples of compact irreducible
$G_2$-manifolds (i.e.\ with holonomy exactly $G_2$) by resolving
singularities of quotients of the flat \mbox{7-torus} by appropriately chosen
finite groups. Later a different type of construction of irreducible
$G_2$-manifolds was developed in~\cite{Ko}, using a {\em generalized
connected sum}. The latter method requires a pair of compact K\"ahler
complex threefolds $\oW_1,\oW_2$ with anticanonical $K3$ divisors satisfying
a certain `matching condition'. Examples of $G_2$-manifolds given
in~\cite{Ko} were constructed from pairs of smooth Fano threefolds, say
$V_1,V_2$, by choosing $K3$ surfaces $S_i\in|-K_{V_i}|$ and blowing up each
$V_i$ along a smooth curve $C_i\in \bigl|-K_{V_i}|_{S_i}\bigr|$.

In this paper, we apply the theory of $K3$ surfaces with non-symplectic
involution~\cite{Ni1, Ni2, Ni3} to give another suitable class of
threefolds $\oW$; these are {\em not} obtainable by blowing up Fano
threefolds along curves as above. Using the connected sum method we
construct hundreds of new topological types of compact  7-manifolds with
holonomy $G_2$, not homeomorphic to those published in~\cite{Jo2} or
\cite{Ko}.
The bounds on Betti numbers attained by our examples are
$0 \leq b^2 \leq 24$, $35 \leq b^3 \leq 239$. Most (but not all) of the
$G_2$-manifolds constructed by Joyce in~\cite[figure~12.3]{Jo1} satisfy
$b^2+b^3\equiv 3\mod 4$. This property does not hold for many of our
examples.

It is an interesting problem to find further classes of threefolds $\oW$
suitable for application of the connected sum method~\cite{Ko}. This is
essentially a problem in algebraic geometry, as can be seen from
Proposition~\ref{alg}. One such class of threefolds arises as a direct
generalization of~\cite[\S 6]{Ko} by blowing up an appropriately chosen
finite sequence of curves, rather than one curve, in a Fano threefold
(Example~\ref{mult}). A full classification of appropriate sequences of
curves extending the classification of non-singular Fano
threefolds~\cite{Is,MoMu,MoMu2} is a separate issue to be dealt with elsewhere.
More generally, a promising strategy is to start with threefolds of
negative Kodaira dimension with vanishing cohomology of their structure
sheaves ($H^{i}(\mco_\oW)=0$ for $0<i\leq 3$) and consider their birational
transformations which have smooth anticanonical divisors. We also note that
Corti, Haskins and Pacini~\cite{CHP} found suitable examples of $\oW$
using weakly Fano threefolds.

The paper is organized as follows. We begin in \S\ref{sum} by giving a
brief introduction to the Riemannian holonomy~$G_2$ and a review of the
conditions on a pair of complex threefolds $\oW_1,\oW_2$ and their
anticanonical $K3$ divisors $D_1,D_2$ required for the connected sum
construction. The section also includes new remarks and observations. In
\S\ref{nsk3} we recall some results about $K3$ surfaces with non-symplectic
involution including the important role of sublattices of the K3 lattice.
These results are applied in \S\ref{3fold} where a class of threefolds
$\oW$ is constructed. We then apply in \S\ref{g2mfds} and \S\ref{further} the
method of~\cite{Ko} and show that, under rather general conditions, the
latter threefolds and the blow-ups of Fano threefolds form pairs suitable
for application of the connected sum construction and discuss the resulting
examples of compact manifolds with holonomy~$G_2$. There is also a short
appendix giving possible values of invariants of Fano threefolds and
$K3$ surfaces with non-symplectic involution.

\section{The connected sum construction of compact irreducible
  $G_2$-manifolds}
\label{sum}

For a detailed introduction to $G_2$-structures see~\cite{Jo2}.
A $G_2$-structure on a 7-manifold $M$ may be given by a real differential
3-form $\phi\in\Omega^3(M)$ such that $\phi$ at every point of $M$ is
isomorphic to
\[
\phi_0=d\bfx_{123}+d\bfx_{145}+d\bfx_{167}+d\bfx_{246}
-d\bfx_{257}-d\bfx_{347}-d\bfx_{356}
\]
on Euclidean~$\RE^7$ with standard coordinates $x_k$, where
$d\bfx_{ijk}=dx_i\we dx_j\we dx_k$. The stabilizer of $\phi_0$ in the
action of $GL(7,\RE)$ is the exceptional Lie group $G_2$ and $\phi$ is
sometimes called a {\em $G_2$-form} on~$M$. As $G_2\subset SO(7)$ every
$G_2$-structure on~$M$ determines an orientation, Riemannian
metric $g(\phi)$ and Hodge star $*_\phi$.

If a $G_2$-form $\phi$ satisfies
\begin{equation}\label{tfree}
d\phi=0\text{ and }d{*_\phi}\phi=0,
\end{equation}
then the holonomy of the metric $g(\phi)$ is contained in~$G_2$
\cite[lemma~11.5]{Sa} and $(M,\phi)$ is called a {\em $G_2$-manifold}. If
$M$ is compact and its fundamental group is finite, then the holonomy
representation of $g(\phi)$ is {\em irreducible} and the holonomy group of
$g(\phi)$ is exactly~$G_2$ \cite[proposition~10.2.2]{Jo1}.

In this section, we give a short summary of the connected sum method
in~\cite{Ko} of constructing solutions to~\eqref{tfree} and further discuss
the class of complex threefolds required for this method, as well as some
topology of the resulting $G_2$-manifolds.

A key ingredient in the connected sum construction of irreducible\linebreak
$G_2$-manifolds is the following result producing
asymptotically cylindrical {\em Calabi--Yau} threefolds.

\begin{theorem}[{\cite[\S 2, \S 3]{Ko}}]\label{bCY}
Let $\oW$ be a compact K\"ahler threefold with K\"ahler form $\omega_{\oW}$
and suppose that $D\in|-K_{\oW}|$ is
a $K3$ surface in the anticanonical class, such that the normal bundle
$N_{D/\oW}$ is holomorphically trivial. Suppose further that $\oW$ is
simply-connected and the fundamental group of $W=\oW\setminus D$ is finite.

Then $W$ admits a complete Ricci-flat K\"ahler metric with
holonomy~$SU(3)$ and a nowhere-vanishing holomorphic $(3,0)$-form.
These are exponentially asymptotic to the product cylindrical Ricci-flat
K\"ahler structure on $D_1\times S^1\times \RE_{>0}$ with K\"ahler form
$\kappa_I+dt\we d\theta$ and non-vanishing holomorphic $(3,0)$-form
$(\kappa_J+\sqrt{-1}\,\kappa_K)\we (dt+\sqrt{-1}d\theta)$.
Here $z=\exp(-t-\sqrt{-1}\,\theta)$ extends over~$D$ to give a holomorphic
local coordinate
on~$\oW$ vanishing to order one precisely on~$D$ and $\kappa_I$ is the
Ricci-flat K\"ahler metric on~$D$ in the class $[\omega_{\oW}|_D]$ and
$\kappa_J+\sqrt{-1}\kappa_K$ is a non-vanishing holomorphic
$(2,0)$-form on~$D$.
\end{theorem}

The next result will be needed in \S\ref{3fold} and can be of
independent interest.

\begin{prop}\label{alg}
Every threefold $\oW$ satisfying the hypotheses of Theorem~\ref{bCY} admits
no non-zero holomorphic forms, $h^{1,0}(\oW)=h^{2,0}(\oW)=h^{3,0}(\oW)=0$. In
particular, $\oW$ is necessarily projective-algebraic.
\end{prop}
\begin{proof}
The vanishing of $h^{1,0}$ is clear as $\oW$ is simply-connected and
K\"ahler. Since the anticanonical bundle of~$\oW$ admits holomorphic
sections (other than the zero section) but $K_{\oW}$ is not holomorphically
trivial, we must have $h^{3,0}=0$. Using the vanishing of
$h^{1,0}$ and $h^{3,0}$ and taking the cohomology of the structure sheaf
sequence for $D\subset\oW$,
\[
0 \ra \mco_\oW(K_{\oW}) \ra \mco_\oW \ra \mco_D \ra 0,
\]
we obtain an exact sequence
\[
0 \ra H^2(\mco_\oW) \ra  H^2(\mco_D) \ra H^3(\mco_\oW(K_{\oW}))\ra 0,
\]
where we also noted that $H^2(\mco_\oW(K_{\oW})) \cong H^1(\mco_\oW) =0$
by Serre duality. Similarly, $\dim H^3(\mco_\oW(K_{\oW}))=1$, so we must
have $H^{2,0}(\oW) =\break  H^2(\mco_\oW) = 0$. The last claim follows
by Kodaira's embedding theorem.
\end{proof}
In what follows, we always assume, as we may, a Ricci-flat K\"ahler metric
and a holomorphic $(3,0)$ form on $\oW$ given by Theorem~\ref{bCY} are
rescaled so that the corresponding triple of 2-forms on a $K3$ surface
$D$ satisfies $\kappa_I^2=\kappa_J^2=\kappa_K^2$.

Let $(W_1,D_1)$ and $(W_2,D_2)$ be two pairs satisfying the assertion of
Theorem~\ref{bCY}. The following definition is crucial for the connected
sum construction of compact irreducible $G_2$-manifolds. It relates to some
standard results about $K3$ surfaces; one excellent reference
is~\cite[chapter~VIII]{BHPV}. Our notation corresponds to the
property that a Ricci-flat K\"ahler structure on a $K3$ surface $D$ is
equivalent to a hyper-K\"ahler structure: there is a triple of complex
structures $I,J,K=IJ=-JI$ such that for each of these the respective
$\kappa_I$, $\kappa_J$, $\kappa_K$ is a K\"ahler form of a Ricci-flat
K\"ahler metric and the two remaining forms span $(H^{2,0}\oplus H^{0,2})(D)$.

\begin{defi}\label{matching}
We say that a pair of Ricci-flat K\"ahler $K3$ surfaces\linebreak
$(D_1,\kappa'_I,\kappa'_J+\sqrt{-1}\kappa'_K)$,
$(D_2,\kappa''_I,\kappa''_J+\sqrt{-1}\kappa''_K)$
satisfies the {\em matching condition} if there is an isometry of
cohomology lattices $h:H^2(D_2,\Z)\to H^2(D_1,\Z)$, so that the
$\RE$-linear extension of $h$ satisfies
\begin{equation}\label{rota}
h: \; [\kappa''_I]\mapsto[\kappa'_J],\quad[\kappa''_J]\mapsto[\kappa'_I],
\quad[\kappa''_K]\mapsto[-\kappa'_K].
\end{equation}
\end{defi}

A map $h$ satisfying the above condition is in fact a pull-back $h=f^*$, for
a well-defined isometry of Riemannian 4-manifolds
$$
f:D_1\to D_2
$$
satisfying
$f\!: \; \kappa''_I\mapsto\kappa'_J,\;\;\kappa''_J\mapsto\kappa'_I,\;\;
\kappa''_K\mapsto -\kappa'_K$ \ \cite[proposition~4.20]{Ko}.
(The map $f$ is sometimes called a `hyper-K\"ahler rotation'; $f$ is
{\em not} in general a holomorphic map with respect to complex structures
on $D_j$ induced by the embeddings $D_j\to\oW_j$, $j=1,2$.)

Each of the two 7-manifolds $W_j\times S^1$ has an asymptotically
cylindrical $G_2$-structure $\phi_j=\omega_j\we d\theta_j+\im\Omega_j$
satisfying~\eqref{tfree} induced by the
asymptotically cylindrical Ricci-flat K\"ahler structure
$\omega_j,\Omega_j$ on~$W_i$. Here $\theta_j$ denotes the standard
coordinate on the~$S^1$ factor. A generalized connected sum of
$W_1\times S^1$ and $W_2\times S^1$ is a compact 7-manifold
\begin{equation}\label{7mfd}
M=(W_1\times S^1)\cup_F(W_2\times S^1),
\end{equation}
defined by truncating the ends of $S^1\times W_j$
at $t=T+1$ to obtain compact manifolds with boundaries and identifying
collar neighbourhoods $D_j\times S^1\times S^1\times [T,T+1]$ of the
boundaries. This latter identification uses a map $F$ of the form
\begin{equation}\label{collar}
\begin{split}
F:D_1\times S^1\times S^1\times [T,T+1]&\to
D_2\times S^1\times S^1\times [T,T+1],
\\
(y,\theta_1,\theta_2,T+t)&\mapsto (f(y),\theta_2,\theta_1,T+1-t) .
\end{split}
\end{equation}
with $f$ as defined above. The map $F$ identifies the $S^1$ factor in
each $S^1\times W_j$ with a circle around the divisor $D_{3-j}$ in
$\oW_{3-j}$, for $j=1,2$, resulting in a finite fundamental group of~$M$
(more precisely, see~(\ref{invariants}a) below).

Furthermore, it is possible to smoothly cut off the $G_2$-forms $\phi_j$
to their asymptotic models on
$D_1\times S^1\times S^1\times [T,T+1]\subset S^1\times W_j$ obtaining
closed $G_2$-forms $\tilde\phi_j$ satisfying
$$
\tilde\phi_j|_{D_j\times S^1\times S^1\times [T,T+1]} =
\kappa_I^{j}\wedge d\theta_j + \kappa_J^{j}\wedge d\theta_{3-j}
+ \kappa_K^{j}\wedge dt + d\theta_{3-j}\wedge d\theta_j\wedge dt.
$$
It follows from~\eqref{collar} that $F^*\phi_2=\phi_1$ on the collar
neighbourhood and so $\tilde\phi_1$ and $\tilde\phi_2$ together give a
well-defined 1-parameter family of closed $G_2$-forms $\phi(T)$ on~$M$.
The corresponding metric $g(\phi(T))$ on $M$ then has diameter asymptotic
to $2T$ when $T\to\infty$. The form $\phi(T)$ is not co-closed but
$\|d*_{\phi(T)}\phi(T)\|\to 0$ as $T\to\infty$. For
every large~$T$, the form $\phi(T)$ can be perturbed into a solution
of~\eqref{tfree} on~$M$.

\begin{theorem}[{\cite[\S 5]{Ko}}]\label{g2}
Suppose that each of $\oW_1,D_1$ and $\oW_2,D_2$ satisfies the hypotheses
of Theorem~\ref{bCY} and the $K3$ surfaces $D_j\in|-K_{\oW_j}|$ satisfy the
matching condition.
Then the compact 7-manifold $M$ defined in~\ref{7mfd} has finite
fundamental group and admits a $G_2$-structure $\phi$
satisfying~\eqref{tfree}, and hence a metric $g(\phi)$ with holonomy equal
to~$G_2$.
\end{theorem}

We can identify some topological invariants of the $G_2$-manifold $M$.
The key to understanding the topology of~$M$ is a non-trivial identification
between $K3$ surfaces $D_1,D_2$ and their second cohomology, given by
Definition~\ref{matching}. The following set-up will be used several times
in the paper. For each $j=1,2$, an embedding $D_j\to\oW_j$ induces a
homomorphism
\begin{equation}\label{iota}
\iota_j:H^2(W_j,\RE)\to H^2(D_j,\RE).
\end{equation}
Let
$$
X=\iota_1(H^2(W_1,\RE))\cap f^*\iota_2(H^2(W_2,\RE)),
$$
and denote
\begin{equation}\label{dn}
d_j=d(\oW_j)=\dim\ker\iota_j,\qquad n=\dim X.
\end{equation}
The restriction of the intersection form on $H^2(D_1,\RE)$ to
$\iota_1(H^2(W_1,\RE))$ is non-degenerate with positive index $1$ because
$W_1$ is K\"ahler, $h^{2,0}(W_1)=0$ and $H^{1,1}(D_1)\cap H^2(D_1,\RE)$ is
non-degenerate has positive index~1. Similarly, considering a subspace
$\iota_2(H^2(W_2,\RE))$ in $H^2(D_2,\RE)$ and noting that
$f^*=h:H^2(D_2,\RE)\to H^2(D_1,\RE)$ preserves the intersection form, we find
that the induced form on $f^*\iota_2(H^2(W_2,\RE))$ is non-degenerate with
positive index~$1$.

A class in $\iota_2(H^2(W_2,\RE))$ with positive square is obtained by
restricting a K\"ahler form on $W_2$; in the notation of Theorem~\ref{bCY}
and Definition~\ref{matching} this class may be taken to be~$\kappa''_I$.
By~\eqref{rota}, $f^*(\kappa''_I)=\kappa_J$ is orthogonal to $H^{1,1}(D_1)$,
thus orthogonal to~$X$. Similarly, $\kappa'_I\in\iota_1(H^2(W_1,\RE))$ is
orthogonal to~$X$. Therefore, $X$ is negative-definite and there are uniquely
defined subspaces $X_j$ so that
\begin{equation}\label{compl}
\iota_1(H^2(W_1,\RE))=X\oplus X_1
\quad\text{ and }\quad
f^*\iota_2(H^2(W_2,\RE)=X\oplus X_2
\end{equation}
are orthogonal direct sum decompositions, with respect to the intersection
form on~$D_1$.

\begin{theorem}\label{topology}
\begin{subequations}\label{invariants}
Let $\oW_j,D_j$ and $M$ be as in Theorem~\ref{g2} with $d_j,n,X_j$ as
defined in \eqref{dn} and~\eqref{compl}. Then
\begin{gather}
\pi_1(M)=\pi_1(W_1)\times\pi_1(W_2),
\\
b^2(M)=n+d_1+d_2.
\end{gather}
Suppose further that $b^2(W_1)-d_1+b^2(W_2)-d_2\le 22$ and that
$X_1$ is orthogonal to $X_2$ with respect to the intersection form on
$H^2(D_1,\RE)$. Then
\begin{equation}
b^3(M) = b^3(\oW_1) + b^3(\oW_2) + b^2(M) - 2n +23.
\end{equation}
\end{subequations}
\end{theorem}

\begin{remark}
A consequence of \eqref{invariants} is that the quantity
$$
b^2(M)+b^3(M) = b^3(\oW_1) + b^3(\oW_2) + 2d_1 + 2d_2 +23
$$
depends only on the threefolds $\oW_j$ but {\em not} on the possible
ambiguity of satisfying the matching condition for $K3$ divisors $D_j$.
\end{remark}

\begin{proof}[Proof of Theorem~\ref{topology}]
This is a direct generalization of \cite[theorems 4.28 and 8.57]{Ko}
and is proved by the same arguments based on application of van Kampen
theorem and Mayer--Vietoris exact sequences. The only difference is that
the computation of Betti numbers in {\it op.cit.} are carried out for a
class of threefolds $\oW_j$ ($j=1,2$) satisfying $d(\oW_j)=0$, but in
general $d(\oW_j)$ need not vanish.

As the normal bundle of $D_j$ is trivial, $D_j$ has a tubular neighbourhood
$U_j$ in $W_j$ so that $U_j\subset W_J$ is homotopy equivalent to
$D\times S^1$. The Mayer--Vietoris for
\begin{equation}\label{MV.W}
\oW_j=W_j\cup U_j
\end{equation}
combined with the K\"unneth formula for $D\times S^1$ and the known
topological invariants of a $K3$ surface $D_j$ yields
\begin{equation}\label{BW}
b^2(\oW_j)=b^2(W_j)+1 \quad\text{ and }\quad
b^3(W_j)=b^3(\oW_j)+22-b^2(W_j)+d_j.
\end{equation}
Here we also used the exact sequence for the absolute $H^*(W_J)$ and relative
de Rham cohomology $H^*_c(W_j)$ to obtain $b^4(W_j)=b^2_c(W_j)=1+d_j$.

The computation of $b^2(M)$ is straightforward from the Mayer--Vietoris
for~\eqref{7mfd}. In order to find $b^3(M)$ we shall determine the
dimension of the kernel of
$$
\tau_j+\iota_{3-j}:H^3(W_j,\RE)\oplus H^2(W_{3-j},\RE)\to H^2(D_j,\RE),
$$
where $\iota_j:H^2(W_j,\RE)\to H^2(D_j,\RE)$ denotes the pull-back via the
inclusion $D_J\subset W_j$ and $\tau:H^3(W_j,\RE)\to H^2(D_j,\RE)$ is a
composition of a similar pull-back and the natural isomorphism
$H^3(D_j\times S^1,\RE)\cong H^2(D_j,\RE)$. From the Mayer--Vietoris
for~\eqref{MV.W} we find $\dim\ker\tau_j=b^3(\oW_j)$.

\begin{lemma}\label{ort}
For each $j=1,2$, there is an orthogonal direct sum decomposition
$$
H^2(D_j,\RE)=\iota_j(H^2(W_j,\RE))\oplus\tau_j(H^3(W_j,\RE))
$$
with respect to the intersection form on $H^2(D_j,\RE)$.
\end{lemma}

Assume Lemma~\ref{ort} for the moment and then $X_j$, for each $j$, must be
a subspace of the image of $\tau_{3-j}$, as $X_{3-j}$ is orthogonal to
$X\oplus X_j$ by the hypothesis.
It follows that $\rk(\tau_j+\iota_{3-j})=n+\rk\tau_j=n+b^3(W_j)-b^3(\oW_j)$
using once again the Mayer--Vietoris for~\eqref{MV.W}.
Now (\ref{invariants}c) follows by rank-nullity for $\tau_j+\iota_{3-j}$,
the K\"unneth formula for $W_j\times S^1$ and the exactness of the
Mayer--Vietoris for~\eqref{7mfd}.

To prove Lemma~\ref{ort} we use the exactness of the sequence for
$H^*_c(W_j,\RE)$ and $H^*(W_j,\RE)$ in the term $H^3(D_j\times S^1,\RE)$
and the Poincar\'e duality. The class $[\omega]\in H^2(D_j,\RE)$ is in the
image of $\tau_j$ if and only if $[\tilde\omega\wedge\omega\wedge
d\theta_j\wedge d\rho(t)]=0$ in $H^6_c(W_j)\cong\RE$, for every class
$[\tilde\omega]\in H^2(W_j,\RE)$. Here $\rho(t)$ denotes a cut-off function
of the cylindrical end coordinate $t$ on $W_j$ (as defined in
Theorem~\ref{bCY}) and, as before, $\theta_j$ is a standard coordinate
on~$S^1$.

The class $[\tilde\omega\wedge\omega\wedge d\theta_j\wedge d\rho(t)]$ is
trivial in $H^6_c(W_j)$ if and only if $[\tilde\omega\wedge\omega]=0$ in
$H^4(D_j,\RE)$ (here we used $H^5(W,\RE)=0$).
This completes the proof of Lemma~\ref{ort} and Theorem~\ref{topology}.
\end{proof}

Suitable examples of $\oW$ are found in~\cite[\S 6]{Ko} by starting with a
(smooth) Fano threefold $V$ with $S$ a smooth anticanonical divisor on~$V$
(then $S$ must be a $K3$ surface). Suppose that $C\in |-K_V|_S|$ is a connected
smooth curve. Blowing up $V$ along $C$, defines a threefold $\oW$ with
$D\subset\oW$ the proper transform of~$S$, so that $D$ is anticanonical
divisor on $\oW$ with holomorphically trivial normal bundle. Thus $\oW$ and
$D$ satisfy the hypotheses of Theorem~\ref{bCY}. Also $d(\oW)=0$ holds in
this case, by the Lefschetz hyperplane theorem.

One instance of non-vanishing $d(\oW)$ arises via the following
generalization of~\cite[\S 6]{Ko}.  Rather than blowing up along one
curve~$C$, we now blow up $V$ along a \emph{finite sequence} of curves.
Suppose that there is a divisor
\[
c=c_1 + \ldots + c_m \in |-K_V|_S|,
\]
where all $c_i$ are connected smooth curves. We do not assume that $c_i$
are distinct. First let $\sigma_1:V_1\to V$ be the blow up $V$ along $c_1$
and consider the proper transform $S_1$ of~$S$. Since $c_1$ lies in $S$,
the $S_1$ is still an anticanonical divisor
on~$V_1$ and $S_1$ is biholomorphic to~$S$ via the restriction of
$\varsigma_1=\sigma_1|_{S_1}$. We use $\varsigma_1$ to identify each $c_i$
with a curve in~$S_1$ still denoted by $c_i$. Then
$$
-K_{V_1}|_{S_1} =
\sigma_1^*(-K_{V})|_{S_1} - c_1 =
\varsigma_1^*(-K_{V}|_{S}) - c_1 =
\varsigma_i^*(-K_{V}|_{S} -c_1) =
c_2 + \ldots + c_m
$$
(cf.~\cite[pp.~604--608]{GH}). Also $V_1$ and $V_1\setminus S_1$ are
simply-connected (cf.~\cite[corollary~6.43]{Ko}).

Next we blow up $c_2$ in~$V_1$ obtaining a threefold $V_2$ with
$-K_{V_2}|_{S_2}=c_3 + \ldots + c_m$. Proceed recursively in this manner,
blowing up $c_i$ in $V_{i-1}$ to obtain $V_i$, up to $i=m$. Denote by
$\oW=V_m$ the result of this sequence of $m$ blow-ups and let
$D=S_m\subset\oW$ be, respectively, the result of the sequence of proper
transforms of~$S$. We find that $D\in|-K_{\oW}|$ and $-K_{\oW}|_{D}=0$,
that is, the normal bundle of $D$ is holomorphically trivial, and
$\pi_1(\oW)=\pi_1(W)=1$, thus $\oW$ and $D$ satisfy the hypotheses of
Theorem~\ref{bCY}. Furthermore,
\[d(\oW)= m - \rk \langle c_1, \cdots, c_m \rangle,\]
where we used $\langle c_1, \cdots, c_m \rangle$ to denote a (free abelian)
subgroup of $H^2(S,\Z)$ generated by $c_i$'s. Thus $d(\oW) > 0$ whenever the
images of $c_1, \cdots, c_m$ are linearly dependent in $H^2(S,\RE)$.

\begin{example}\label{mult}
Let $V=\CP^3$. Then $S$ can be taken to be any quartic $K3$ surface.
Let $c=h_1 + \cdots + h_4$, where $h_i$'s are smooth hyperplane
sections of $S$. We can construct a threefold $\oW$ as above; in the
present case $m=4$ and $d(\oW) = 3$. We shall return to this example
in~\S\ref{seq}.
\end{example}

\section{Non-symplectic involutions on $K3$ surfaces}
\label{nsk3}

A holomorphic involution $\rho$ of a $K3$ surface $S$ is called
{\em non-symplectic} when $\rho^*(\omega) = -\omega$ for each
$\omega \in H^{2,0}(S)$. Non-symplectic involutions of $K3$ surfaces were
studied and completely classified by V.V. Nikulin (\cite{Ni1, Ni2, Ni3},
see also \cite[\S 6.3]{AlNi}). In this section, we recall,
in summary, some results that will be needed in this paper.

The intersection form of $S$ makes the integral cohomology $H^2(S,\Z)$ into
a {\em lattice}, i.e.\ a free abelian group of finite rank endowed with
an integral symmetric bilinear form. Geometric properties of~$S$ depend
strongly on the data of its second cohomology. The arithmetic of the lattice
$H^2(S,\Z)$ and its sublattices plays a key role in Nikulin's results. We
begin with a brief account of the required foundational concepts of
lattice theory and refer the reader to~\cite[\S~I.2]{BHPV} and \cite{Do1}
for more details.

We shall use $x\cdot y \in\Z$ to denote the bilinear form evaluated on two
lattice elements $x,y$ and sometimes write $x^2$ for the associated
quadratic form. The lattice is called {\em even} if $x^2$ is even for every
lattice vector~$X$. The notions of signature, non-degenerate lattice,
orthogonal complement and lattice isometries (isomorphisms) are understood
relative to the bilinear form of this lattice. By the direct sum of
lattices we shall always mean the orthogonal direct sum. The lattice is
called unimodular when the matrix of its bilinear form in some basis has
determinant $\pm 1$.

A sublattice $A$ of a lattice $N$ is {\em primitive} if the quotient group
$N/A$ has no torsion or, equivalently, any basis of $A$ extends to a basis
of~$N$. In particular, $x\in N$ is a primitive element if $x$ generates a
primitive sublattice, i.e.\ not divisible (in $N$) by an integer $>1$.
If a lattice is isometric to a primitive sublattice of~$N$ then
it is said to be {\em primitively embedded} in~$N$.

For any non-degenerate lattice $N$ the dual abelian group
$N^*=\Hom_{\Z}(N,\Z)$ naturally contains $N$ as a subgroup, via the
inclusion map $x\in N\to (x(y):=y\cdot x) \in N^*$. The quotient
$N^*/N$ is called the {\em discriminant group} of~$N$; this group is
finite, of order $|\det (e_i\cdot e_j)|$, for any basis $\{e_i\}$ of~$N$.
We denote by $l(N)$ the minimal number of generators of $N^*/N$. This
quantity satisfies
\begin{equation}\label{le.rank}
l(N)\le \rk N
\end{equation}
and will be important several times later.
We shall sometimes need the induced `dual' quadratic form on
$N^*$ with values in~$\Q$ and the induced quadratic form $N^*/N\to
\Q/(2\Z)$, defined in the usual way (here $N$ is assumed even).

It is well-known that for each $K3$ surface $S$ the lattice $H^2(S,\Z)$ is
non-degenerate, with signature $(3,19)$, unimodular, even and (hence)
isomorphic to $L=2(-E_8)\oplus 3 H$. Here $E_8$ denotes the unique (up to
isometry) even, unimodular and positive definite lattice of rank~8 and
$H=\bigl(\begin{smallmatrix} 0&1\\ 1&0 \end{smallmatrix}\bigr)$ is the
standard hyperbolic plane lattice. We shall refer to $L$ as the {\em $K3$
lattice} and denote by $O(L)$ its group of isometries. Any lattices
considered in this paper will be isomorphic to sublattices of~$L$.

Denote by $L^\rho$ the set of all classes in $H^2(S,\Z)$ fixed by $\rho^*$.
It is clear that $L^\rho$ is a primitive sublattice of $H^2(S,\Z)$ and,
since $h^{2,0}(S)=1$, it also follows that $L^\rho$ is a sublattice
of the Picard lattice $H^{1,1}(S)\cap H^2(S,\Z)$. If $h$ is a class of some
K\"ahler form on~$S$, then $h+\rho^*(h)$ is again a K\"ahler class (as
$\rho$ is holomorphic) which is also invariant under~$\rho^*$. It is then
easy to see that $L^\rho$ must contain an (integral) K\"ahler class, thus
every $K3$ surface $S$ with non-symplectic involution $\rho$ is algebraic.

We shall interchangeably think of $\rho^*$ as an involution of~$L$ (an
isometry of order 2 in $O(L)$) and $L^\rho$ as a sublattice of~$L$. The
sublattice $L^\rho$ is called the {\em invariant lattice} (of $\rho$) and
is a primitive non-degenerate (by the Hodge index theorem) sublattice
of~$L$ with signature $(1,t_-)$. Let $r=1+t_-$ denote its rank.

The orthogonal complement of $L^\rho$ is the `anti-invariant' sublattice,
i.e.\ precisely the $(-1)$-eigenspace of $\rho^*$. As $L$ is unimodular,
thus $L^*\cong L$, the ($\Z$-linear) orthogonal projections $L\to L_+^*$ and
$L\to L_-^*$ induce isomorphisms of groups
\begin{equation}\label{2iso}
L_+^*/L_+ \cong L/(L_+\oplus L_-)\cong L_-^*/L_- .
\end{equation}
which are $\rho$-equivariant. Then every element of the group
in~\eqref{2iso} must have order~2; a group with this property is called
{\em 2-elementary}.

In particular, for each involution $\rho$, the discriminant group of the
invariant lattice $(L^\rho)^*/ L^\rho$ is 2-elementary, isomorphic to
$(\Z_2)^a$. Clearly $l(L^\rho)=a$ in this case.

Finally, define
\begin{equation}\label{delta}
\delta(L^\rho) =
\left \{ \begin{array}{ll}
  0 & \text{ if } t^2 \in \mbz \text{ for each } t \in {L^\rho}^*, \\
  1 & \text{ otherwise,}
\end{array} \right.
\end{equation}
using the quadratic form on ${L^\rho}^*$ with values in~$\Q$ induced by the
bilinear form of~$L$. Remark that the definition of $\delta$ makes sense
for each non-degenerate lattice.

\begin{prop}[{\cite[\S 4]{Ni3}}]\label{clas}
The triple $(r,a,\delta)$ defined above forms a complete system of invariants
which determine the involution $\rho^*$ and $L^\rho$ uniquely up to isometries
of~$L$.
\end{prop}
We shall denote by $L{(r, a, \delta)}$ the invariant lattice determined by
$(r,a,\delta)$.

Conversely, suppose that $N$ is a non-degenerate lattice with signature
$(1,t_-)$, primitively embedded in the $K3$ lattice $L$ and such that $N^*/N$
is 2-elementary. Then there is an involution $\rho_N$ on $L$ with the
invariant lattice precisely~$N$. Then $\rho_N$ acts as $-1$ on the
orthogonal complement $N^\bot$ of~$N$ in~$L$ and the $\RE$-linear
extension of $N^\bot$ in $L\otimes\RE$ contains a positive 2-plane, $P$ say.
By the global Torelli theorem and the surjectivity of the period
map for $K3$ surfaces~\cite[chapter~VIII]{BHPV}, there is a well-defined
(non-singular) $K3$ surface $S$ with non-symplectic involution $\rho$, so
that $\rho^*=\rho_N$ and $H^2(S,\RE)\cap (H^{2,0}\oplus H^{0,2})(S)$ is
identified with~$P$ via (the $\RE$-linear extension of) some lattice
isomorphism $H^2(S,\Z)\to L$.

That $S$ is non-singular can be seen by noting that $L_{Pic}=P^\bot \cap L$
corresponds to the Picard lattice of~$S$ and the orthogonal complement of $N$
in $L_{Pic}$ (relative to the bilinear form of~$L$) cannot contain elements
$x$ with $x^2=-2$. Indeed, otherwise the Riemann--Roch for the divisor $x$
on $S$ gives that $x$ or $-x$ is effective, but $\rho^*(x)=-x$, which is absurd.

With the above understood, deformations families of $K3$ surfaces with
non-symplectic involution are identified with isomorphism classes of
primitive 2-elementary $(1,t_-)$-sub\-lattices of~$L$ with invariants
$(r,a,\delta)$.
All the triples $(r,a,\delta)$ that occur for some $L^\rho$ were
classified in \cite{Ni3}. There are 75 possibilities; the range of values is
$1 \leq r \leq 20$ and $0 \leq a \leq 11$ (see Appendix).
By \eqref{le.rank}, any possible triple $(r,a,\delta)$ satisfies
\begin{equation}\label{ineq}
r-a \ge 0.
\end{equation}

For a lattice $N$, it will be convenient to consider,
following~\cite{Do2}, the concept of {\em ample $N$-polarized} $K3$
surface defined by the conditions that its Picard lattice contains a
sublattice isomorphic to~$N$ and this sublattice contains a class
represented by some K\"ahler form (sometimes also called an ample
class). Ample polarized $K3$ surfaces may be thought of as a refinement of the
concept of algebraic $K3$ surfaces of degree $2n-2$ ($n\ge 2$) in $\CP^n$;
in this latter case, the lattice $N$ is generated by one
element (a K\"ahler form) of positive square $2n-2$.
The moduli space for ample $N$-polarized $K3$ surface is a quasiprojective
variety of complex dimension $20-\rk N$ and is irreducible whenever
$N^\bot\subset L$ contains a copy of the hyperbolic plane lattice~$H$,
cf.~\cite[proposition~5.6]{Do2}.

As noted above, every $K3$ surface $S$ with non-symplectic involution with
invariants $(r,a,\delta)$ has an integral ample class in the invariant
lattice. Thus $S$ is ample $L{(r,a,\delta)}$-polarized.
Denote by $\KK{(r, a, \delta)}:=\KK_{L{(r, a, \delta)}}$ the moduli space
of $L{(r, a, \delta)}$ lattice polarized $K3$ surfaces and by
$\KK'{(r, a, \delta)}$ a subspace of $\KK{(r, a, \delta)}$, whose
elements are $K3$ surfaces with non-symplectic involution with
the invariant lattice isomorphic to ${L}{(r,a,\delta)}$.

\begin{prop}[{\cite[\S 6.3]{AlNi} or \cite[\S 4]{Ni3}}]\label{klem}
\mbox{}
\begin{enumerate}
\item
$\KK'{(r,a,\delta)}$ is connected and  Zariski open in~$\KK{(r,a,\delta)}$.
\item
for $(S,\rho) \in \KK'{(r,a,\delta)}$, the fixed locus of $\rho$ is
\begin{itemize}
\item $S^\rho = \emptyset$ if ${(r, a, \delta)} = (10, 10, 0)$;
\item $S^\rho$ is a disjoint union of two elliptic curves if
 ${(r, a,\delta)} = (10, 8, 0)$;
\item otherwise
\[
S^\rho =c_g + e_1 + \cdots + e_k,\quad\text{ all disjoint,}
\]
where $g = (22 - r - a)/2$, $k = (r - a)/2$,  $c_g$ is a curve
of genus $g$ and $e_i \cong \CP^1$.
\end{itemize}
\end{enumerate}
\end{prop}
\begin{remark}
The invariant $\delta$ also has a geometric interpretation related to
the above proposition: it is shown in~\cite{AlNi} that $\delta(L^\rho)$ is
zero if and only if the fixed locus $S^\rho$ of $\rho$ on~$S$ represents a
class divisible by~2 in $H^2(S,\mbz)$.
\end{remark}

\section{Threefolds $W$ from non-symplectic involutions}
\label{3fold}

Let $\psi:\CP^1 \ra \CP^1$ be any holomorphic involution fixing two distinct
points and define
$Z = (S \times \CP^1) /(\rho, \psi)$. Then the singular locus of $Z$ is a
product of smooth curves and ordinary double points, resulting from the
fixed locus of $\rho$. Let $\overline W \ra Z$ be the blow-up along the
singular locus of $Z$. It is elementary to check that $\overline W$ is
smooth.  Choose a point $p \in \CP^1$ such that $p\neq \psi(p) $. Let $D'$
be the image of $S\times \{p\}$ in $Z$ and $D$ be the inverse image of $D'$
in $\overline W$. Then $D$ is isomorphic to $S$ and it is an anticanonical
divisor of $\overline W$ with the normal bundle $N_{D/\oW}$ holomorphically
trivial. One can obtain $\overline W$ differently. Let ${q_1, q_2}$ be the
fixed points of $\psi$, $S^\rho = c_1 + \cdots + c_l$ and $Q_{ij}$ be $ c_i
\times \{q_j\} \subset S \times \CP^1$. Let $\widetilde W$ be the blow-up
of $S\times\CP^1$ along the curves $Q_{ij}$'s. Then the involution on
$S\times\CP^1$ induces an involution on $\widetilde W$, whose fixed locus
are the exceptional divisors over $Q_{ij}$. The quotient of~
$\widetilde W$ by the involution is isomorphic to $\overline W$. This may
be summarized by
$$
\begin{CD}
\widetilde{W} 	@>>> \oW\\
@VVV 		@VVV\\
S\times\CP^1 	@>>> Z
\end{CD}
$$
which is an instance of resolution of singularities diagram, for $Z$.

It follows from Lemma~\ref{klem}(a) that for any fixed
$(r,a,\delta)$ the threefolds $\oW$ constructed from
$S\in\KK'{(r,a,\delta)}$ are deformations of each other.  The following
fact will be used later:

\begin{prop} \label{ampleh}
Suppose that $\kappa \in L{(r,a,\delta)} \otimes \mbr \subset H^2(D, \mbr)$
is a K\"ahler class. Then $\kappa$ is a restriction of a K\"ahler
class in $H^2(\overline W, \mbr)$.
\end{prop}
\begin{proof}
Let $E_1, \cdots, E_r$ be the exceptional divisors of the blow-up
$g:\widetilde W \ra S \times \mbc P^1$. Then $g^*(\kappa, \mco_{\mbc
P^1}(1))-\varepsilon (E_1 + \cdots + E_r)$ is a K\"ahler class of
$\widetilde W$ for sufficiently small positive~$\varepsilon$.
Since this class is invariant under the induced involution on
$\widetilde W$, there is a K\"ahler class $\hat \kappa \in
H^2(\overline W, \mbr)$ whose pullback via the quotient map
$\widetilde W \ra \overline W$ is $g^*(\kappa, \mco_{\mbc
P^1}(1))-\varepsilon (E_1 + \cdots + E_r)$. Clearly the restriction
of $\hat \kappa$ to $D$ is $\kappa$.
\end{proof}

Our next lemma ensures that Theorem~\ref{bCY} can be applied to $\oW$
and~$D$ defined above.
\begin{lemma}\label{W}
$\overline W$ and $W:=\overline W \setminus D$ are simply connected if
$\rho$ is not fixed-point-free, i.e.\ if $(r, a, \delta) \neq (10,10,0)$.
\end{lemma}
\begin{proof}
Let $\widehat W \ra \overline W$ be the universal covering. The threefold
$\widetilde W$ is simply-connected, therefore the quotient map
$\widetilde W\to\oW$ lifts to a holomorphic map $\widetilde{W} \ra \widehat{W}$
so there is a commutative diagram:
\[
\xymatrix{
    &   \widehat W \ar[d]\\
\widetilde W \ar[r]\ar[ur] & \overline W .}
\]
If $\rho$ has a fixed point, then $\widetilde W \ra \widehat W$ cannot be a map
of degree one, so it is a map of degree two. Then $\widehat W \ra \overline W$
is of degree one and $\widehat W$ is necessarily isomorphic to $\overline W$,
thus $\overline W$ is simply-connected.

There is a natural projection: $(S\times\CP^1)/(\rho, \psi) \ra S /\rho$.
Let $\nu$ be the composition of $\overline W \ra Z$ and $(S\times\CP^1)/(\rho,
\psi) \ra S /\rho$. Let $x \in S/\rho$ be a point in the branch locus of
the map $S \ra S /\rho$. Then $\nu^{-1}(x)$ is a union of three smooth
rational curves, one of which (denoted by $l$) crosses $D$ transversely at
a single point and the other two are disjoint from $D$, resulting from the
blow-up . Since $\overline W$ and $D$ are simply connected, $\pi_1(W)$ is
generated by a loop around $D$. We can assume that the loop is contained in
$l^* = l - D$. Since the loop can be contracted to a point in $l^*$,
$W$ is simply-connected.
\end{proof}
Since we are only interested in the case that $\oW$ is simply connected, we
will assume that  $(r, a, \delta) \neq (10, 10, 0)$ in the rest of this
paper.

We require Betti numbers of $\oW$ for application of Theorem~\ref{topology}
later.  In light of Proposition~\ref{alg} it suffices to find $h^{1,1}$ and
$h^{1,2}$.
Let $\widetilde X \ra X$ be the blow-up of smooth variety along a smooth
subvariety $Z$.  The following formula for the topological Euler numbers is
well-known:
\begin{equation} \label{euler1}
e(\widetilde X) = e(X)+e(E)-e(Z),
\end{equation}
where $E$ is the exceptional divisor.
Let $X \ra Y$ be a double covering over smooth variety $Y$, branched along
smooth subvariety $B$. The following formula is also well-known:
\begin{equation} \label{euler2}
e(X)=2 e(Y) - e(B).
\end{equation}
\begin{prop}\label{hodge}
A threefold~$\oW$ constructed from a $K3$ surface in $\KK'(r, a, \delta)$,
\linebreak
$(r, a, \delta)\neq (10,10,0)$, as defined above, satisfies:
\begin{enumerate}
\item
 $h^{1,1}(\oW)=b^2(\oW)=3+2r-a$  and $h^{1,2}(\oW)=\frac12\, b^3(\oW)=22-r-a$.
\item $\rk\;\iota =r$, where  $\iota:H^2(\overline W,\RE) \ra H^2(D,\RE)$
  is the restriction map (cf.~\eqref{iota}).
\end{enumerate}
\end{prop}
\begin{proof}
\begin{enumerate}
\item \label{ccc} the Hodge number
$h^{1,1}(\overline W)$ can be obtained by counting the number
of independent classes in $H^{1,1}(\widetilde W)$ that are invariant under
the involution. Let $S^\rho = c_1 + \cdots + c_l$. We have, by standard
results about blow-ups (e.g.~\cite[p.~605]{GH}),
\[
H^2(\widetilde W,\CX)\cong H^2(S \times \CP^1,\CX)\oplus\CX^{2l}
=H^2(S,\CX)\oplus\CX\oplus\CX^{2l},
\]
using also the K\"unneth formula. The last two terms correspond to
$H^{1,1}(\CP^1)$ and the exceptional divisors; all of these classes are
fixed by the involution. By the hypothesis, exactly $r$ independent classes
in $H^{1,1}(S)$ are invariant under the involution, so we obtain
\[
h^{1,1}(\oW) = r + 1 + 2 l.
\]
By Equation (\ref{euler1})), the topological Euler number of $\widetilde W$ is
\[e(\widetilde W) = 24 \cdot 2 + 2 \sum_{i} e(c_i)\]
On the other hand, we have
\[2 e(\oW) -\sum_{i,j} e (Q_{ij})= e(\widetilde W)\]
by Equation (\ref{euler2}).
So we have
\begin{align*}
e(\oW) &=\frac{1}{2} \left (\sum_{i,j} e (Q_{ij}) + e(\widetilde W) \right )\\
       &=2 \sum_i e(c_i)  + 24  +  \sum_{i} e(c_i)\\
       &=24 + 3 \sum_i e(c_i).
\end{align*}
Finally we have
\[h^{1,2}(\oW)=1+h^{1,1}(\oW) - \frac{1}{2} e(\oW) =-10+r+2l - \frac{3}{2}   \sum_i e(c_i).\]
\item since $h^{2,0}(\oW)=0$, we have
$\rk\, (H^2(\overline W,\RE) \ra H^2(D,\RE))=\rk_{\CX}  \varrho$,
where we denoted by $\varrho:H^{1,1}(\overline W) \ra H^{1,1}(D)$
the restriction map.
As seen in the proof of (a), $H^{1,1}(\oW) = A\oplus\CX^{1+2l}$, where
$A\cong L(r,a,\delta)\otimes\CX$ corresponds to the classes invariant under
the involution. It is easy to see that $\dim \varrho(A)=r$ and
$\varrho(\CX^{1+2l})=0$. Thus $\dim \im \varrho = r$ as required.
\end{enumerate}
\end{proof}
As a corollary we also compute the quantity defined in Theorem~\ref{topology}
\begin{equation}\label{wd}
\begin{split}
d(\overline W) &= \dim\ker (H^2(\overline W, \RE) \ra H^2(D, \RE)) - 1\\
&=\dim\ker (H^2(W,\RE) \ra H^2(D,\RE))\;
=\;      2+r-a
\end{split}
\end{equation}
In particular, the kernel of $H^2(\overline W, \mbz) \ra H^2(D, \mbz)$
contains the classes corresponding to exceptional divisors of the blow-up
$\oW\to Z$.

We find, taking account of~\eqref{ineq}, that the threefolds $\oW$
constructed in this section have $d(\oW)\ge 2$. In light of the remark at
the end of~\S\ref{sum}, these threefolds are {\em not} homeomorphic to any
threefold obtainable by blowing up a curve in a Fano threefold as
in~\cite[\S 6]{Ko}. They also cannot be homeomorphic to the threefold in
the Example~\ref{mult} as the latter has $d=3$ but according to the
classification of $K3$ surfaces with non-symplectic involution
$2+r-a$ is always even.

\section{Compact irreducible $G_2$-manifolds}
\label{g2mfds}

For each $K3$ surface $S$ with non-symplectic involution having fixed
points, we constructed in the previous section a threefold $\oW$ with
$D\in|-K_{\oW}|$ biholomorphic to~$S$ and such that the normal bundle of
$D$ is holomorphically trivial.  For a fixed isomorphism class of invariant
lattice $L(r,a,\delta)$, the resulting pairs $(\oW,D)$ are deformations of
each other. In particular, it follows from the construction and
Proposition~\ref{klem}(a) that the isomorphism classes of $K3$ surfaces $D$
arising in these pairs form a {\em Zariski open} (in particular, open and
dense) subset $\KK'(L(r,a,\delta))\subset\KK(L(r,a,\delta))$ of the space of
ample $L(r,a,\delta)$-polarized $K3$ surfaces. For brevity, we shall sometimes
refer to such $(\oW,D)$ and the respective polarizing lattice $L(r,a,\delta)$
to be of {\em non-symplectic type}.

By Lemma~\ref{W} both $\oW$ and $W$ are simply-connected and by
Proposition~\ref{ampleh} there exists a K\"ahler form $\omega_{\oW}$ on
$\oW$ such that $[\kappa_I]=[\omega_{\oW}|_D]$ in $H^2(D,\RE)$. Thus
Theorem~\ref{bCY} applies to give the following.

\begin{prop}\label{prepNS}
Let $(\oW,D)$ be of non-symplectic type and
$(r, a, \delta) \neq (10,10,0)$. Suppose that
$\kappa_J+\sqrt{-1}\kappa_K$ is a non-vanishing holomorphic $(2,0)$-form
on~$D$ and $\kappa_I$ is the K\"ahler form of some Ricci-flat K\"ahler
metric on~$D$. Assume $\kappa_I^2=\kappa_J^2=\kappa_K^2$.

Then $W=\oW\setminus D$ has an asymptotically cylindrical Ricci-flat K\"ahler
metric and holomorphic $(3,0)$-form modelled on the product structure
on $D\times S^1\times \RE_{>0}$ determined by
$\kappa_I,\kappa_J+\sqrt{-1}\kappa_K$ as defined in Theorem~\ref{bCY}.
\end{prop}

Another class of pairs $(\oW,D)$ satisfying the hypotheses of
Theorem~\ref{bCY} was given in~\cite[\S 6]{Ko}. Recall from \S\ref{sum}
that we start with a Fano threefold~$V$ and consider a $K3$ divisor $S\in
|-K_V|$. Blowing-up a connected smooth curve in $V$ representing the
self-intersection cycle $S\cdot S$ (multiplication in the Chow ring) yields
a threefold $\oW$ and the proper transform $D$ of~$S$, so that the pair
$(\oW,D)$ satisfies the hypotheses of Theorem~\ref{bCY}. Furthermore, $D$ and
$S$ are isomorphic {\em as K\"ahler surfaces}, with K\"ahler forms induced
from respectively $\oW$ and $V$ and with appropriate choice of K\"ahler form
on~$\oW$. By varying $V$ in its algebraic family and varying $S$ in the
anticanonical class we have an ambiguity to choose $D$ in a Zariski open
subset $\KK'(L(V))\subset\KK(L(V))$ of the space of ample $L(V)$-polarized
$K3$ surfaces \cite[theorem~6.45]{Ko}. Here the sublattice
\begin{equation}\label{fano.sublat}
L(V)=\iota_V^*(H^2(V,\Z))\subset H^2(D,\Z)
\end{equation}
is induced by the embedding $\iota_V:D\to V$ and $\rk L(V)=b^2(V)$, as in
the present case, $\iota_V^*$ is injective by the Lefschetz hyperplane
theorem. The resulting pairs $(\oW,D)$ are deformations of each other; we
stress that the lattice $L(V)$ is a topological invariant, independent of
the choice of $(\oW,D)$ in its deformation class. We shall refer to such
$(\oW,D)$ and the respective polarizing lattice $L(V)$ to be
of {\em Fano type} (though $\oW$ is never a Fano threefold).

The resulting class of asymptotically cylindrical Calabi--Yau threefolds
may be described as follows.

\begin{prop}[{cf.~\cite[pp.~146--147]{Ko}}]\label{prepF}
Let $(\oW,D)$ be of Fano type, obtained from Fano threefold~$V$. Let
$\kappa_J+\sqrt{-1}\kappa_K$ be a non-vanishing holomorphic $(2,0)$-form
on~$D$. Suppose that $\kappa_I$ is the K\"ahler form of a Ricci-flat K\"ahler
metric on~$D$, such that the cohomology class $[\kappa_I]$ is induced from
some K\"ahler form on~$V$ via the proper transform $D\to S\subset V$.
Assume $\kappa_I^2=\kappa_J^2=\kappa_K^2$.

Then $W=\oW\setminus D$ has an asymptotically cylindrical Ricci-flat K\"ahler
metric and holomorphic $(3,0)$-form modelled on the product structure
on $D\times S^1\times \RE_{>0}$ determined by
$\kappa_I,\kappa_J+\sqrt{-1}\kappa_K$ as defined in Theorem~\ref{bCY}.
\end{prop}

\begin{remark}
For notational convenience, we restrict the term `Fano type' to
threefolds obtainable by blow-up of {\em one} connected curve. Thus the
threefold in Example~\ref{mult} is not of Fano type (though it is a very
natural generalization); we discuss this threefold separately.
\end{remark}

In this section, we shall slightly abuse the notation by writing
$\KK'(L_j)$ for both the non-symplectic and Fano types when, respectively,
$L_j=L(V)$ and $L_j=L(r,a,\delta)$. However, as we shall
only be using the property that $\KK'(L_j)$ is some Zariski open
subset of $\KK(L_j)$ this should not cause confusion.

If by deforming two pairs $(\oW_1,D_1)$, $(\oW_2,D_2)$ satisfying
Theorem~\ref{bCY} and rescaling, if necessary, the respective K\"ahler
metrics and holomorphic forms we achieve the matching condition for
$D_1$ and $D_2$ (Definition~\ref{matching}), then by Theorem~\ref{g2},
the connected sum construction described in \S\ref{sum} yields a compact
Riemannian 7-manifold $M$ with holonomy $G_2$. We shall sometimes
abbreviate this by saying that a compact (irreducible) $G_2$-manifold $M$
is {\em constructed from $(\oW_1,D_1)$ and $(\oW_2,D_2)$} (or simply from
$\oW_1$ and $\oW_2$, the existence of appropriate $D_1$ and $D_2$ implied).

Here is the main theorem that will be applied to construct new examples of
compact 7-manifolds with holonomy~$G_2$ in this paper.

\begin{theorem}\label{keythm}
For $j=1,2$, let $(\oW_j,D_j)$ satisfy the hypotheses of Theorem~\ref{bCY}
with $L_j$ the respective sublattice of the $K3$ lattice.
Suppose that, for each $j$, the $K3$ surfaces arising from deformations of
the pair $(\oW_j,D_j)$ define a Zariski open subset in the moduli space
$\KK(L_j)$ of $L_j$-polarized $K3$ surfaces. 
Suppose further that at least one of the following two conditions holds:
\begin{enumerate}
\item
the lattice $L_1 \oplus L_2$ can be embedded in the $K3$ lattice so that
the restriction to each $L_j$ is a primitive embedding or
\item
each $(\oW_j,D_j)$ is of Fano type or of non-symplectic type
as defined above and, 
if $L_j$ is of non-symplectic type with invariants $(r_j,a_j,\delta_j)$,
then either $(r_j,a_j,\delta_j)=(10,8,0)$ or $r_j\le 5$.
\end{enumerate}
Then a compact irreducible $G_2$-manifold $M$ can be constructed from
$(\oW_1,D_1)$ and $(\oW_2,D_2)$, in the sense defined above.
Furthermore, $M$ may be taken to satisfy
$b^2(M)=d_1+d_2$, \ i.e.\ with $n=0$, in the notation of~\eqref{dn}.
\end{theorem}

Theorem~\ref{keythm} is a modification of~\cite[theorem~6.44]{Ko}. More
explicitly, the latter result asserts that a compact irreducible
$G_2$-manifold can be constructed from $(\oW_j,D_j)$ when both pairs are of
Fano type, which is a situation allowed by clause (b) of
Theorem~\ref{keythm}. Note that this includes the Fano threefold discovered
in~\cite{MoMu2} and not appearing in the list of~\cite{MoMu} or \cite{IsPr}
cited in~\cite{Ko}.
On the other hand, the argument in~\cite[\S 6]{Ko} is carried out in such a
way that it only uses the arithmetic of sublattices $L_j\subset L$ and the
result that, for $(W_j,D_j)$ of Fano type, $\KK'(L_j)$ is a connected
Zariski open subset of $\KK(L_j)$ (Proposition~\ref{prepF}). As we already
noted in Proposition~\ref{prepNS} a similar property holds when $(W_j,D_j)$ is
of non-symplectic type.

Unfortunately, there is a mistake in Lemma 6.47 in~\cite{Ko} which is a
part of the argument of Theorem~6.44 in {\it op.cit.} in the case when the
polarizing sublattices have $\rk L_j>1$. We therefore include a complete
proof of Theorem~\ref{keythm} avoiding any dependence on the affected
lemma.

We begin with lattice embedding results of Nikulin~\cite{Ni2} which will
be used in the proof of Theorem~\ref{keythm} and again in~\S\ref{examples}.

\begin{theorem}\label{lattice}
A primitive embedding of an even non-degenerate lattice $N$ of signature
$(t_+,t_-)$ into an even unimodular lattice $E$ of signature $(l_+,l_-)$ exists
provided $t_+\leq l_+$, $t_-\leq l_-$ and at least one of the following
conditions holds:
\begin{enumerate}
\item
$2\rk N \le \rk E$ or
\item
$\rk N + l(N) < \rk E$, where $l(N)$
is the minimum number of generators of the discriminant group $N^*/N$.
\end{enumerate}

If moreover $t_+< l_+$, $t_-< l_-$ and one of the following two
inequalities holds, $\rk N + l(N) \le \rk E - 2$ or $2\rk N \le \rk E - 2$,
then a primitive embedding $N\to E$ is unique, up to isometry of $E$.
\end{theorem}

Clause (a) is gathered from \cite[theorem~1.12.4]{Ni2} and
(b) from \cite[corollary~1.12.3]{Ni2}, see also a simplified version given
by Dolgachev \cite[theorem (1.4.6)]{Do1}. For the uniqueness part,
see~\cite[theorem 14.4.4]{Ni2}. We also note that if an even
unimodular lattice $E$ is indefinite, then $E$ is determined by its
signature $(l_+,l_-)$ uniquely up to isometry.

\begin{remark}
The existence condition (b) in Theorem~\ref{lattice} can be slightly
weakened. For each prime~$p$, let $l(N)_p$ denote the minimal number of
generators of the Sylow $p$-subgroup of $N^*/N$. Then $l(N)=\max_p l(N)_p$.
In the notation of Theorem~\ref{lattice}, the lattice $N$ admits a
primitive embedding into~$E$ whenever $\rk N + l(N)_2 \le \rk E$ and $\rk N
+ l(N)_p < \rk E$ for each odd prime~$p$.
\end{remark}

We also require an instance of primitive embedding not directly obtainable
from Theorem~\ref{lattice}.

\begin{prop}\label{embed}
(i) Let a lattice $L_j$ be of Fano type and $\rk L_j\ge 6$. Then
$\rk L_j\le 10$ and $L_j$ is determined by its rank uniquely up to
isomorphism and admits a primitive embedding in $-E_8(2)\oplus H$.
If $\rk L_j=10$, then $L_j$ is isomorphic to the non-symplectic
type lattice $L(10,8,0)\cong -E_8(2)\oplus H$.

(ii) Suppose that the polarizing lattice $L_j$ of a $K3$ divisor
$D_j\subset\oW_j$ has $\rk L_j\le 5$.
Then  $L_j$ admits a primitive embedding in $-E_8\oplus H$ and,
if $\rk L_j\le 4$, also a primitive embedding in $-E_8(2)\oplus 2H$.

Here $-E_8(2)$ denotes the negative-definite lattice of rank~8 obtained by
multiplying the bilinear form of $-E_8$ by~$2$.
\end{prop}
\begin{proof}
(i) Every Fano threefold $V$ satisfies $b^2(V)\le 10$. If
$b^2(V)\ge 6$, then $V$ is biholomorphic to a product $\CP^1\times
S_{11-b^2(V)}$ \cite{MoMu,MoMu2}. Here $S_d$ denotes a del Pezzo surface of
degree $d$ biholomorphic to a blow-up of $\CP^2$ in
$9-d$ generic points.

The following properties of the lattice $L(S_d\times\CP^1)$ are obtained by
more-or-less standard methods. The intersection form on $L_j=L(V)$ may be
computed as a cup product $x\smallsmile y\smallsmile c_1(V)$ in
$H^{\text{even}}(V,\Z)$. The generators of $H^2(S_{11-b^2(V)}\times\CP^1)$
are Poincar\'e dual to the 4-cycles of $\CP^2\times\mathrm{pt}$,
$\ell\times\CP^1$ and $E_i\times\CP^1$, where $\ell$ corresponds to a
projective line in $\CP^2$ and $E_i$ are the exceptional divisors on~$S_d$;
the Chern class $c_1$ of the threefold corresponds to $(2,3,-1,\ldots,-1)$
in this basis. Straightforward computations yield
\begin{equation}\label{delpezzo}
L(S_d\times\CP^1)=
\begin{pmatrix}
0 & 3 & \mathbf{1}^t\\
3 & 2 & \mathbf{0}^t\\
\mathbf{1} & \mathbf{0} & I_m(-2)
\end{pmatrix}.
\end{equation}
Here $m=9-d$ and we used $\mathbf{1}$ and $\mathbf{0}$ to denote a column
vector in $\Z^m$ with all components equal, respectively, to~$1$ and to~$0$
and $I_m(-2)$ is an $m\times m$ diagonal matrix whose all diagonal entries
are~$-2$.

We find that $L(S_d\times\CP^1)$, for each $1\le d\le 8$, is a primitive
sublattice of $L(S_1\times\CP^1)$, thus it suffices to consider
\eqref{delpezzo} with $d=1$, $m=8$. The order of the discriminant 
group of $L(S_1\times\CP^1)$ is 256, obtained by taking the determinant
of~\eqref{delpezzo}. The rows of~\eqref{delpezzo} give a basis of the image
of $L(S_1\times\CP^1)$ in the dual $\Z$-module $L(S_1\times\CP^1)^*$; the
corresponding reduced basis (in the sense of Lovasz, computed e.g.\ by
the `LLL algorithm' for integer bilinear forms~\cite{LLL}) has an upper
triangular form $\{e_1,\,\sum_{i=2}^{10} e_i,\, 2e_j\,|\, j=3,\ldots,10\}$
(here $e_j$ is the standard basis of~$\Z^{10}$). It is then not difficult to
check that the discriminant group of $L(S_1\times\CP^1)$ is a 2-elementary
group $(\Z_2)^8$. Therefore, by Proposition~\ref{clas}, $L(S_1\times\CP^1)$
must be isomorphic to a non-symplectic type lattice $L(10,8,0)$ as $\delta=0$
is the only possibility in this case according to the classification of
non-symplectic involutions on $K3$ surfaces, see Appendix. (It can also be
checked directly that the dual quadratic form $q$ on $L(S_1\times\CP^1)^*$
takes only integer values.) We obtain
$$
L(S_1\times\CP^1)\cong L(10,8,0)\cong -E_8(2)\oplus H
$$
and the maps
\begin{equation}\label{e8}
\iota_\pm:e\in -E_8(2) \to (e,\pm e)\in 2(-E_8),
\end{equation}
with each choice of sign, define primitive embeddings. These induce primitive
embeddings $L(S_1\times\CP^1)\to 2(-E_8)\oplus H$ acting as identity on the
$H$ component.

(ii) If $L_j=L(V)$ with $b^2(V)\le 5$ or $L_j=L(r,a,\delta)$ with $r\le 5$,
then $L_j$ has a primitive embedding in $-E_8\oplus H$, by
Theorem~\ref{lattice}.

Now suppose that $\rk L_j\le 4$. To obtain a primitive
embedding of $L_j$ in a non-unimodular lattice $-E_8(2)\oplus 2H$ we follow
the method of~\cite[proposition 1.15.1]{Ni2}.

The lattice $-E_8(2)\oplus 2H$ is a primitive sublattice of a
unimodular $(-2)E_8\oplus 2H$, with the orthogonal complement isomorphic
to~$-E_8(2)$. Suppose that $S$ is a lattice containing $L_j$ as a primitive
sublattice such that the orthogonal complement of $L_j$ in~$S$ is
isomorphic to $-E_8(2)$. Thus $L_j$ and $-E_8(2)$ are primitive sublattices
of~$S$ and are orthogonal complements of each other and the quotient group
$S/\bigl(L_j\oplus -E_8(2)\bigr)$ is finite. The discriminant group
$S^*/S$ is a subgroup of $(L_j^*/L_j)\oplus (-E_8(2)^*/-E_8(2))$. The
latter quotient space has a natural 
quadratic form $q=q_L\oplus q_E$ with values in $\Q/2\Z$, where $q_L$ and
$q_E$ are induced by the quadratic forms of, respectively, $L_j$ and
$-E_8(2)$. It is not difficult to check that $q|_{S^*/S}=0$. On the other
hand, since the sublattices $L_j$ and $-E_8(2)$ are primitive in~$S$ the
projections $p_L: S/\bigl(L_j\oplus -E_8(2)\bigr)\to L_j^*/L_j$ 
and  $p_E: S/\bigl(L_j\oplus -E_8(2)\bigr)\to (-E_8(2)^*/-E_8(2))$ are
injective group homomorphisms. Thus $q\circ p_L=-q\circ p_E$.
Conversely, the ambiguity of choosing a lattice $S$ as above is equivalent
to a choice of an isomorphism $\gamma:H_L\to H_E$ of subgroups
$H_L\subset L_j^*/L_j$ and $H_E\subset (-E_8(2))^*/-E_8(2)$ satisfying
$q_E\circ\gamma=-q_L$. The discriminant group of $S$ is then given by
$$
S^*/S\cong
\bigl(\bigl((L_j^*\oplus -E_8(2)^*)/H\bigr)/(L_j\oplus -E_8(2))\bigr)/H,
$$
where $H$ is the graph of $\gamma$ (cf.~\cite[propositions 1.4.1
and~1.5.1]{Ni2}).

The minimal number of generators of $L_j^*/L_j$ satisfies $l(L_j)\le\rk
L_j\le 4$. Let $G_L$ be the Sylow $2$-subgroup of $L_j^*/L_j$. Then
$G_L$ has a subgroup $H_L\cong (\Z_2)^{a_L}$ with $a_L=l(L_j)_2\le 4$. We
can obtain an injective homomorphism $\gamma:H_L\to (-E_8(2))^*/(-E_8(2))
\cong (\Z_2)^8$ satisfying $q_E\circ\gamma=-q_L$, where $H_E=\gamma(H_L)$.
Here we used $l(L_j)_2\le l(-E_8(2))/2$ and straightforward
calculation based on the `Jordan decomposition' for quadratic forms over
the field $\Z_2$; we omit the details but see e.g.~\cite[\S~9.2]{AlNi}.
The resulting lattice $S$ has $l(S)_2\le 9$, $l(S)_p\le 4$ for each
prime $p\neq 2$ and admits a primitive embedding $\iota_S:S\to
(-2)E_8\oplus 2H$ by the remark after Theorem~\ref{lattice}.

The restriction of $\iota_S$ to $-E_8(2)$ defines a primitive embedding
$-E_8(2)\to 2(-E_8)\oplus 2H$ which is unique by
Theorem~\ref{lattice}. By composing with an isometry of
$2(-E_8)\oplus 2H$ we may assume that $\iota_S$ maps $-E_8(2)$ onto the
orthogonal complement of $-E_8(2)\oplus 2H$ in $(-2)E_8\oplus 2H$. Then the
restriction of $\iota_S$ to $L_j$ gives a primitive embedding
$L_j\to -E_8(2)\oplus 2H$.
This completes the proof of Proposition~\ref{embed}.
\end{proof}

\begin{proof}[Proof of Theorem~\ref{keythm}.]
We use the argument of~\cite[pp.~149--150]{Ko}
adapting it to the present situation, in combination with
Proposition~\ref{embed}, to give the proof of Theorem~\ref{keythm} in
the case of hypothesis~(b). A simple modification of this argument
will also prove the result in the case~(a).

By the hypotheses, Theorem~\ref{bCY} holds for each $(W_j,D_j)$. Therefore,
as discussed above, it suffices to show that the matching condition for the
triples of cohomology classes
$\kappa^{(j)}_I,\kappa^{(j)}_J,\kappa^{(j)}_K\in H^2(D_j,\Z)$
can be achieved by some deformations of $(W_j,D_j)$ (here $j=1$ or $2$
corresponds to $'$ or $''$ in~\eqref{rota} and we dropped square brackets
to avoid excessive notation). Our result will then follow by application of
Theorem~\ref{g2}.

For each of the $K3$ surfaces $D_j$ we have a choice a lattice isometry
(often called a marking) $p_j:H^2(D_j,\Z)\to L=2(-E_8) \oplus 3H$.
Clearly, the ambiguity of choosing $p_j$ is the group of lattice
automorphisms of~$L$. Recall that the K\"ahler $K3$ surfaces $D_j$ are ample
\mbox{$L_j$-polarized}: $L_j\subseteq \iota_j^* H^2(\oW_j,\Z)\subset
H^2(D_j,\Z)$ consists of $(1,1)$-classes. We may assume that the K\"ahler
class $\kappa^{(j)}_I$ is primitive in $H^2(D_j,\Z)$. The restriction of
$p_j$ identifies $L_j$ with a sublattice $p_j(L_j)$ of~$L$. In terms of
`abstract' lattices, this defines a primitive embedding $p_j:L_j\to
L$. Conversely, any two primitive embeddings $L_j\to L$ are related by an
automorphism of~$L$ (we proved a stronger statement in
Proposition~\ref{embed}). Thus every primitive embedding $L_j\to L$ arises from
some marking~$p_j$.

The subspace $P_j=H^2(D_j,\RE)\cap (H^{2,0}\oplus H^{0,2})(D_j)$ is spanned
by $\kappa^{(j)}_J$ and $\kappa^{(j)}_K$ and is identified via~$p_j$ with a
positive-definite real 2-plane in $L\otimes\RE$. This 2-plane is
necessarily orthogonal to the sublattice $p_j(L_j)$.  Recall from the
discussion in \S\ref{nsk3} that the converse also holds, by the global
Torelli theorem and the surjectivity of the period map for $K3$ surfaces.
More precisely, every positive-definite real 2-plane $P$ in $L\otimes\RE$
orthogonal to $p_j(L_j)$ is the image of
$H^2(D_j,\RE)\cap (H^{2,0}\oplus H^{0,2})(D_j)$,
for some ample \mbox{$L_j$-polarized} $K3$ surface in $\KK(L_j)$, provided
that $P^\bot\cap L$, which is the image of the Picard lattice, does not
contain any elements with square $-2$ orthogonal to a K\"ahler class.
(Remember that we use the bilinear form on~$L$ throughout these
considerations.) As we noted earlier, the $K3$ surfaces arising from
deformations of pairs $(\oW_j,D_j)$ form a Zariski open subset
$\KK'(L_j)\subset\KK(L_j)$. Thus we are free to choose $P_j$ in a Zariski
open subset of the Grassmannian of positive 2-planes orthogonal to $p_j(L_j)$.

The task of satisfying the matching condition is now equivalently stated
in terms of sublattices of the $K3$ lattice. We shall sometimes write
$L=(-E_8)_1\oplus(-E_8)_2\oplus H_1\oplus H_2\oplus H_3$ to distinguish
between the copies of $(-E_8)$ and $H$ in the direct sum.

Applying Theorem~\ref{lattice} and Proposition~\ref{embed}, we obtain for
each $L_j$ satisfying the hypothesis (b) a primitive embedding
in $-E_8(2)\oplus H$ when $\rk L_j\ge 6$ or a primitive embedding in
$-E_8\oplus H$ when $\rk L_j\le 5$ or a primitive embedding in
$-E_8(2)\oplus 2H$ when $\rk L_j\le 4$. In each case, the embedding is
induced by some marking $p_j$ of $D_j$, $j=1,2$.

If each polarizing lattice has $\rk L_j\le 5$, then we choose a primitive
embedding $L_j\to (-E_8)_j\oplus H_j$ for each $j=1,2$. If each $L_j$ has
$\rk L_j\ge 6$, then there is a primitive embedding of $L_1$ in
$-E_8(2)_+ \oplus H_1$ and $L_2$ in $-E_8(2)_- \oplus H_2$. Here we used
$-E_8(2)_\pm$ to denote the images $\{(e,\pm e):e\in -E_8\}\subset 2(-E_8)$
of the primitive embedding $\iota_\pm$ defined in~\eqref{e8}. If say $\rk
L_1\ge 6$ and $\rk L_2\le 4$, then $L_1$ is primitively embedded in $-E_8(2)_+
\oplus H_1$ and $L_2$ in $-E_8(2)_- \oplus H_2\oplus H_3$.

Finally, suppose that $\rk L_1\ge 6$ and $\rk L_2=5$. As before, we use
a primitive embedding $L_1\to -E_8(2)_+ \oplus H_1$. According to the
classification results, there are exactly three possibilities for $L_2$ of
Fano type and two of non-symplectic type.

One such $L_2=L(V)$ of Fano type arises from $V=S_6\times\CP^1$; in
this case, $L_2$ has a primitive embedding in $-E_8(2)_- \oplus H_2$.
The other two cases, $V$ is obtained by successively blowing up four
projective lines $\CP^1$ (satisfying appropriate conditions),
respectively, in $\CP^3$ or in a smooth quadric $Q\subset\CP^4$ \cite{MoMu}.
A straightforward, though slightly lengthy, calculation shows that
the discriminant groups of $L(V)$ is $\Z_4\oplus\Z_{11}$, for the blow-up
of $\CP^3$, and $(\Z_2)^3\oplus\Z_7$ for the blow-up of~$Q$. (Both groups
can be determined from the prime factorization of the determinant of
$L_2$, noting also that  $l(L_2)_2 \equiv \rk L_2\mod 2$
cf.~\cite[proposition~11.1.4]{Ni2}.) In the $\Z_4\oplus\Z_{11}$ case, we
have $l(L_2)=1$ and, by the splitting theorem \cite[corollary~1.13.5]{Ni2},
$L_2\cong H\oplus L'$ for some negative-definite 
rank 3 lattice~$L'$ with $l(L')=1$. We can obtain a primitive embedding
$L'\to -E_8(2)\oplus H$ by applying, twice, an argument similar to that of
Proposition~\ref{embed}(ii). First, we choose a rank 4 lattice $S'$
containing $L'$ as a primitive sublattice and such that $(S')^*/S'\cong
(\Z_2)^2\oplus H'$, where $H'$ has no elements of even order and $l(S')=
l(S')_2=2$. Then we choose a lattice $S''$ containing $S'\oplus -E_8(2)$,
such that $S''/(S'\oplus -E_8(2))$ is finite and $l(S')=l(S')_2=6$,
$l(S'')_p\le 2$. A primitive embedding of $L'\to -E_8(2)\oplus H$ is
obtained as restriction of a primitive embedding $S''\to 2(-E_8)\oplus H$
which exists by the remark after Theorem~\ref{lattice}. It follows that
$L'\oplus H$ primitively embeds in $(-E_8(2))_-\oplus H_2\oplus H_3$.

In the case when $L_2^*/L_2$ is $(\Z_2)^3\oplus\Z_7$, an argument similar
to Proposition~\ref{embed}(ii) shows that we can choose a lattice
$S$ containing $L_2\oplus -E_8(2)$ with $S/(L_2\oplus -E_8(2))$ finite and
with $l(S)_2\le 7$, $l(S)_p=1$. Then $S$ embeds in $2(-E_8)\oplus 2H$ by
the remark after Theorem~\ref{lattice}, so $L_2$ primitively embeds in
$(-E_8(2))_-\oplus H_2\oplus H_3$.

If $L_2$ is of non-symplectic type then either $L_2\cong L(5,5,1)\cong
I_4(-2)\oplus\mathbf{1}(2)$ or $L_2\cong L(5,3,1)\cong I_3(-2)\oplus H$
(see Appendix), where, as before, $I_m(-2)$ denotes a diagonal rank~$m$
matrix with entries $-2$. Noting from~\eqref{delpezzo}, with $m=8$, that
$-E_8(-2)\oplus H$ has a primitive sublattice isomorphic to $I_8(-2)$ ,
it is easy to see that each of the $L(5,5,1)$, $L(5,3,1)$ has a primitive
embedding in $-E_8(2)_- \oplus H_2\oplus H_3$.

To sum up, in each situation, we have the following orthogonality property
of the primitive embeddings $p_j:L_j\to L$. Firstly, the images of K\"ahler
classes of $D_1$, $D_2$ satisfy
\begin{equation}\label{compat}
p_1(L_1)\bot p_2(\kappa^{(2)}_I)\qquad\text{and}\qquad
p_2(L_2)\bot p_1(\kappa^{(1)}_I).
\end{equation}
Secondly, consider the sublattices
\begin{subequations}\label{orthog}
\begin{equation}
N_3=\bigl(p_1(L_1)\cup p_2(L_2)\bigr)^\bot,
\qquad
N_j=\bigl(p_j(L_j)\oplus N_3\bigr)^\bot.
\end{equation}
It is easy to see that
\begin{equation}
N_j\subseteq p_{3-j}(L_{3-j}),
\qquad
N_j\oplus N_3\subseteq p_j(L_j)^\bot
\text{ and }
\rk p_j(L_j)^\bot =\rk N_j + \rk N_3.
\end{equation}
\end{subequations}
and that the signatures of $N_j$, $N_3$ are, respectively, $(1,s_j)$, $(1,s)$.
Thus the orthogonal complement of $p_1(L_1)\cup p_2(L_2)$ in $L$ contains
positive vectors.
(We actually have a stronger property $p_1(L_1)\bot p_2(L_2)$, though this
is only required for showing $b^2(M)=d_1+d_2$ and will not be used in the
anywhere else in the proof.)

Choose a positive vector $v\in \bigl(p_1(L_1)\cup p_2(L_2)\bigr)^\bot$
and denote $\kappa_j=p_j(\kappa_I^{(j)})$.
For each $j=1,2$, the positive 2-plane $P_j$ through $v$ and $\kappa_j$
defines a $K3$ surface in $\KK(L_j)$, though possibly not in the required
Zariski open subset $\KK'(L_j)$. In order to ensure that the $K3$ surfaces
are in the latter Zariski open subspaces one more step is necessary.

The family of positive 2-planes through a positive $v\in N_3\otimes\RE$ and a
positive $v_j\in N_j\otimes\RE$ defines a real-analytic subvariety
$\KK(L_j)^{\RE}$ of real dimension $20-\rk L_j$ in the moduli space $\KK(L_j)$
of complex dimension $20-\rk L_j$. The real subvariety $\KK(L_j)^\RE$ is
locally modelled on the product of projective spaces $P^+_j\times P^+$ of real
lines in the positive cones of $N_j\otimes\RE$ and $N_3\otimes\RE$,
On the other hand, the local Torelli theorem implies that the complex
$(20-\rk L_j)$-dimensional moduli space $\KK(L_j)$ is locally biholomorphic
to the period domain $\{\omega\CX\in \mathbb{P}(L_j^\bot\otimes\CX):
(\omega,\omega)=0, (\omega,\bar{\omega})>0\}$.
We deduce by examination of the tangent spaces that
$\KK(L_j)^\RE$ is {\em not} contained in any complex subvariety of
positive codimension in~$\KK(L_j)$. Thus $\KK(L_j)^\RE$ cannot lie
in any Zariski closed subset of $\KK(L_j)$ and has to contain points
defined by `generic' $K3$ surfaces in $\KK'(L_j)$. The complement of the
respective generic period points in $\KK(L_j)^\RE$ is the common zero set
of a finite system of real-analytic functions. By the identity theorem,
this analytic subset cannot contain open neighbourhoods, so the points in
$\KK(L_j)^\RE$ defined by generic $K3$ surfaces form a dense open subset.

In particular, any `small' neighbourhood $U$ of $v_3$ in $N_3\otimes\RE$
admits a dense open $U_j\subseteq U$, so that for each $v\in U_j$ there is
a $\tka_j\in N_j\otimes\RE$ `close' to $\kappa_j$ with the property that
the positive 2-plane through $v$ and $\tka_j$ defines a $K3$ surface in
$\KK'(L_j)$. (Here we used `small' and `close' in the sense of Euclidean
norm on $N_j\otimes\RE$, temporarily ignoring the bilinear form on~$L$.)
We claim that with $v$ chosen in the intersection of two respective dense open
subsets $U_1\cap U_2$ we achieve the matching condition.

More precisely, let $j=1$. The perturbed $K3$ surface $D_1$ in $\KK'(L_1)$
is defined by a small perturbation of the 2-plane
$\RE\tka_2\oplus\RE v$ which is still positive, so $\CX(\tka_2+iv)$ is a
period point and determines the complex structure on~$D_1$.
Replacing the K\"ahler class $\kappa_1$ of~$D_1$ with $\tka_1$ corresponds
to a small deformation $\tka^{(1)}_I$ of the K\"ahler class on $D_1$ (here
we noted $N_j\subseteq p_{3-j}(L_{3-j})$). This lifts to a small
deformation of the K\"ahler form on the compact threefold $\oW_1$ containing
$D_1$ as an anticanonical divisor because the pull-back map of $H^2(W,\Z)\to
L_1$ (composed with the marking) is surjective and the
K\"ahler property is preserved by small deformations.
Note that a K\"ahler metric corresponding
to $\tka_1$ need not come from any projective embedding, even though $D_1$ is
an algebraic surface, as the K\"ahler class $\tka_1$ need not be integral.
The same argument applies {\it mutatis-mutandis} to~$D_2$. We thus obtain
well-defined `generic' ample $L_j$-polarized, K\"ahler $K3$ surfaces 
$D_j\in\KK'(L_j)$ with K\"ahler classes $\tka_j$ and period points
$\tka_{3-j}+(-1)^{3-j} iv$.

The positive vector $v$ is orthogonal to both $p_1(L_1)$ and
$p_2(L_2)$ is such that each positive 2-plane spanned by
$p_j(\tka^{(j)}_I)$ and $v$ defines a $K3$ surface arising
as some deformation of $(\oW_j,D_j)$, by Propositions~\ref{prepNS} or
\ref{prepF}.  The $\RE$-linear extensions $p_j:H^2(D_j,\RE)\to
L\otimes\RE$ of the markings of $D_j$ satisfy
$$
p_1(\tka^{(1)}_I)=p_2(\tka^{(2)}_J),\quad
p_1(\tka^{(1)}_J)=p_2(\tka^{(2)}_I),\quad
p_1(\tka^{(1)}_K)=-p_2(\tka^{(2)}_K)=v
$$
(compare \eqref{rota}), which is the required matching condition between
$K3$ divisors $D_1$ and $D_2$ on some deformations of $\oW_1$ and~$\oW_2$.
This completes the proof of Theorem~\ref{keythm}(b).
\medskip

If instead the condition~(a) holds, then the sublattice $L_1$ is orthogonal to
$L_2$ by the hypothesis and $N_j=L_{3-j}$, in the above notation. Then
\eqref{compat} is automatically true for the given 
embedding of $L_1\oplus L_2$ in the $K3$ lattice. The restriction of this
embedding to $L_j$ is a primitive embedding $L_j\to L$ which, as discussed
earlier, is induced by some well-defined marking~$p_j:H^2(D_j,\Z)\to L$. We
then obtain a positive vector $v\in H_3$  and proceed with the rest of the
argument as in the case~(b).
\end{proof}

\begin{remarks}
The Betti numbers of a $G_2$-manifold
obtained via (b) are not in general uniquely determined by those of
$W_1,W_2$, but also depend on the choice of matching (see
Theorem~\ref{topology}). In the proof of Theorem~\ref{keythm}(b) the matching
was chosen with $p_1(L_1)\cap p_2(L_2)=\{0\}$ as this simplifies some
details of the lattice arithmetic. However, the argument remains valid when 
$p_1(L_1)$ and $p_2(L_2)$ intersect non-trivially, provided only that the
orthogonality conditions \eqref{compat} and \eqref{orthog} hold.

When none of the conditions (a),(b) in Theorem~\ref{keythm} is
satisfied, it is sometimes still possible to achieve matching pairs
$D_1,D_2$. One such example is identified in~\cite{KN}, where both
$(\oW_j,D_j)$ are of non-symplectic type with the invariant lattices
$L_j=L(10,8,0)$, but the $\RE$-linear extensions of their images in $L$
intersect in a subspace of dimension~4.
\end{remarks}

\subsection{Examples}
\label{examples}

We construct examples of compact $G_2$-manifolds and compute their Betti
numbers using Theorem~\ref{keythm}(a). If (a) holds, then the hypotheses of
Theorem~\ref{topology} are satisfied with $n=0$ and the formulae for the
Betti numbers of resulting $G_2$-manifold become
\begin{subequations}\label{betti}
\begin{align}
b^2(M)&=d_1+d_2,\\
b^3(M)&= b^3(\oW_1) + b^3(\oW_2) + d_1+d_2 + 23.
\end{align}
\end{subequations}
If $\oW_j$ is of Fano type, then $d_j=0$ by application of Lefschetz
hyperplane theorem and
\begin{equation}\label{b3F}
b^3(\oW_j)= b^3(V_j)+(-K_{V_j}^3)+2
\end{equation}
where $V_j$ is the initial Fano threefold (see~\cite[p.~154]{Ko}). For
$\oW_j$ of non-symplectic type, we computed the quantities $d_j$ and
$b^3(\oW_j)$ in Proposition~\ref{hodge}.

If an embedding $L_1\oplus L_2\to 2(-E_8) \oplus 3H$
is primitive, then so is the restriction to each~$L_j$.
The following is a corollary of the embedding Theorem~\ref{lattice}.
\begin{prop}\label{ND}
A primitive embedding of an even non-degenerate lattice $N$ of signature
$(t_+,t_-)$ into the $K3$ lattice $2(-E_8) \oplus 3H$ exists provided
$t_-\leq 19$, $t_+\leq 3$ and one of the next two inequalities holds,
$\rk N\le 11$ or $l(N) + \rk N < 22$, where $l(N)$ is the minimal number
of generators of the discriminant group $N^*/N$.
\end{prop}

Noting that
\begin{enumerate}
\item $l(N) \leq \rk N$,
\item $l(N_1\oplus N_2) = l(N_1) + l(N_2)$ and
\item $l(L{(r, a, \delta)})=a$,
\end{enumerate}
and taking account of Lemma~\ref{W}, Proposition~\ref{hodge} and
formulae~\eqref{betti}, we obtain our next result, a sufficient condition
for application of Theorem~\ref{keythm}(a).

\begin{theorem} \label{emb}
The lattice $L_1 \oplus L_2$, defined in Theorem~\ref{keythm}, can be
primitively embedded in the $K3$ lattice $2(-E_8) \oplus 3H$ in the
following cases ($r_j=\rk L_j$):
\begin{enumerate}
\item
$\oW_1$ and $\oW_2$ are of Fano type, respectively $V_1$ and $V_2$, and
\[
r_1 + r_2 \le 11.
\]
Then the resulting $G_2$-manifold $M$ has
\[
b^2(M)=0\; \text{ and }\;
b^3(M) = b^3(V_1) - K_{V_1}^3 + b^3(V_2) - K_{V_2}^3 + 27,
\]
where $V_1$,$V_2$ are the initial Fano threefolds.
\item
$\oW_1$ is of Fano type $V_1$ and $\oW_2$ is of non-symplectic type
$(r_2,a_2,\delta_2)$ and
\[
r_1 + r_2 \le 11\; \text{ or }\; 2 r_1 + r_2 + a_2 < 22.
\]
Then the resulting $G_2$-manifold $M$ has
\[b^2(M) = 2 + r_2 -a_2\;  \text{ and }\;
b^3(M) = b^3(V_1) -K_{V_1}^3 - r_2 - 3 a_2 + 71.
\]
\item
$\oW_1, \oW_2$ are of non-symplectic types $(r_1,a_1,\delta_1)$,
$(r_2,a_2,\delta_2)$ and
\begin{gather*}
r_1 + r_2 \le 11\; \text{ or }\; r_1 + r_2 + a_1+a_2 < 22.
\end{gather*}
Then the resulting $G_2$-manifold $M$ has
\begin{equation}\label{c}
b^2(M) = 4 + r_1 + r_2 - (a_1 + a_2)\;  \text{ and }\;
b^3(M) = 115 - (r_1 + r_2) - 3 (a_1 + a_2).
\end{equation}
\end{enumerate}
All of the above $G_2$-manifolds are simply-connected, hence irreducible.
\end{theorem}
\begin{remark}
The case (a) with $r_1=1$ is a part of~\cite[theorem~8.57]{Ko}.
\end{remark}

The Betti numbers of $G_2$-manifolds given by Theorem~\ref{emb} are readily
computed from the classification of Fano threefolds in \cite{Is,MoMu}
(or~\cite[chapter~12]{IsPr}) and \cite{MoMu2} and the classification of
$K3$ surfaces with non-symplectic involution in \cite{Ni1}
or~\cite{AlNi} (see Appendix). These Betti numbers satisfy
$0 \le b^2(M)\le 24$ and 
$53 \le b^3(M)\le 239$. The values of $b^2(M)$ are {\em even}, as all the
values of $r-a$ realized by the invariant lattices $L(r,a,\delta)$ are even.
The values of $b^3(M)$ are {\em odd} as also both $b^3(V)=2h^{1,2}(V)$ and
$-K_V^3$ are even for each Fano threefold~$V$. This parity of the Betti
numbers of~$M$ is a consequence of the particular way to achieve a matching
of the $K3$ divisors chosen in Theorem~\ref{keythm}(a) and does not hold in
general for irreducible \mbox{$G_2$-manifolds} obtainable by Theorem~\ref{g2}.
We give examples with other parities of $b^2(M),b^3(M)$ in~\S\ref{further}.

\begin{table}[h]
\caption{Betti numbers of $G_2$-manifolds from Theorem~\ref{emb}}
\begin{center}\begin{tabular}{| c | c | c |}
\hline
\multirow{2}{*}{\phantom{m}$b^2$\phantom{m}}&
number of&
\multirow{2}{*}{values of $b^3$}\\
& pairs $(b^2,b^3)$&\\
\hline
\multirow{2}{*}{0}&
\multirow{2}{*}{70}&
 55, \ldots, 147,\\
&& 151, \ldots, 189,195,197,239\\
\hline
\multirow{2}{*}{2}&
\multirow{2}{*}{49}&
 53, \ldots, 131,\\
&& 135, 141, 145, 149, 153, 157, 161, 165, 169, 173\\
\hline
\multirow{2}{*}{4}&
\multirow{2}{*}{47}& 59, \ldots, 133,\\
&& 141, 147, 151, 155, 159, 163, 167, 171, 175\\
\hline
\multirow{2}{*}{6}&
\multirow{2}{*}{39}&    61, \ldots, 123,\\
&&   143, 147, 149, 153, 157, 161, 165\\
\hline
\multirow{2}{*}{8}&
\multirow{2}{*}{36}&  63, \ldots, 121,\\
&&   145, 151, 153, 155, 159, 163\\
\hline
\multirow{2}{*}{10}&
\multirow{2}{*}{36}&   65, \ldots, 123,\\
&&  147, 153, 157, 159, 161, 165\\
\hline
\multirow{2}{*}{12}&
\multirow{2}{*}{36}& 67, \ldots, 125,\\
&&   149, 155, 159, 163, 165, 167\\
\hline
\multirow{2}{*}{14}&
\multirow{2}{*}{23}&
69, \ldots, 93,\\
&& 97, 99, 101, 105, 107, 111, 113, 115, 151, 157\\
\hline
{16}&
{14}&
71, 73, 75, 79, 83, 85, 87, 89, 91, 95, 99, 103, 113, 153\\
\hline
{18}&
{14}&
73, 75, 77, 81, 85, 87, 89, 91, 93, 97, 101, 105, 115, 155\\
\hline
{20}&
{13}&
75, 77, 79, 83, 87, 89, 91, 93, 99, 103, 107, 117, 157\\
\hline
{22}&
{1}&
93\\
\hline
{24}&
{1}&
95\\
\hline
\end{tabular}\end{center}
\label{tbl}
\end{table}
The $G_2$-manifolds with $b^2(M)=0$ given by Theorem~\ref{emb} have
$H^2(M,\Z)=0$ and are precisely those coming from clause (a).
Many, but not all of these are obtainable from~\cite[theorem~8.57]{Ko};
neither overlap with \cite{Jo1} which only has one example with $b^2=0$ (it
has $b^3=215$). The $G_2$-manifolds obtainable via clause (b) have
$b^2=2,4,\ldots,20$ and those via clause (c) have $b^2=4,6,\ldots,24$.
We list all these pairs $(b^2,b^3)$ in Table~\ref{tbl}. There and later in the
paper the dots denote either consecutive odd integers or consecutive even
integers, according to the context.
Some of these pairs also occur in $G_2$-manifolds constructed by Joyce
in~\cite{Jo1}, but most of the examples in~\cite{Jo1} (239 out of 252) satisfy
an additional relation $b^2 + b^3 \equiv 3 \mod 4$. Theorem~\ref{emb} gives
many new topological types which do not satisfy the latter relation between
Betti numbers.

Altogether we obtained over 8,000 pairs of threefolds $\oW_1,\oW_2$ with
matching $K3$ divisors and 379 distinct pairs of values of Betti numbers
$(b^2,b^3)$ realized by the $G_2$-manifolds from Theorem~\ref{emb}.
Of these, 324 pairs $(b^2,b^3)$ are different from those of $G_2$-manifolds
in \cite[figure~11.1]{Jo2} (note that the latter Figure includes
5 examples from~\cite{Ko}, also obtainable by Theorems~\ref{emb} or
\ref{keythm}).

\section{Further examples of compact $G_2$-manifolds}
\label{further}
In this section, we give three more families of examples of $G_2$-manifolds,
realizing 84 new values of $(b^2,b^3)$, obtained by application of
Theorem~\ref{keythm}, but not coming from Theorem~\ref{emb}.

\subsection{From mirror pairs of $K3$ surfaces to compact $G_2$-manifolds}
\label{mirror}
We shall assume $r+a \neq 22$ and $(r,a,\delta) \neq (10,10,0)$ and let
$L_1=L{(r, a, \delta)}$, $L_2=L{(20-r, a, \delta)}$ throughout this
subsection. There are 36 such pairs of lattices. These pairs of invariant
lattices are called {\em mirror pairs} and were successfully exploited
independently by C.~Voisin~\cite{Vo} and C.~Borcea~\cite{Bo} in the
construction of mirror pairs of Calabi--Yau threefolds. Mirror families of
$K3$ surfaces were also studied by Dolgachev~\cite{Do2}.

Borcea and Voisin in particular noticed that there exists a primitive
embedding of $L_1,L_2$ into the $K3$ lattice $L$ such that
$L_j^{\bot} = L_{3-j} \oplus H$ for
$j=1,2$. Thus $L_1\oplus L_2$ can be embedded in
$2(-E_8)\oplus 2H$ is satisfied and
a compact irreducible $G_2$-manifold~$M$ can be constructed from
the respective pair of threefolds $\oW_j$ of non-symplectic type
$(r, a, \delta)$ and $(20-r, a, \delta)$ and Theorem~\ref{emb}(c)
applies to give
$$
b^2(M)=24-2a,\qquad b^3(M)=95-6a.
$$
In this way, we obtain 11 distinct pairs of Betti numbers $(b^2, b^3)$ of
$G_2$-manifolds:
\begin{equation}\label{mir}
\begin{split}
(4, &35),
(6, 41),
(8, 47),
(10, 53),
(12, 59),
(14, 65),\\
&(16, 71),
(18, 77),
(20, 83),
(22, 89),
(24, 95),
\end{split}
\end{equation}
satisfying a linear relation
\[
b^3 = 3 b^2 + 23.
\]
which implies $b^2+b^3\equiv 3\mod 4$ satisfied by many $G_2$-manifolds
in~\cite{Jo2}. The last three pairs in~\eqref{mir} do not appear in the
list~\cite[Fig.~11.1]{Jo2}. Two of these are in Table~1; the example
$(22, 89)$, to our knowledge, is new.

\subsection{Using the threefold from Example~\ref{mult}}
\label{seq}
Let $\oW_1=\oW$ be the obtained by a consecutive blow-up of four plane quartic
curves $c_i$ in $\CP^3$ as defined in Example~\ref{mult}. For each blow-up the
exceptional divisor $E$ is a fibre bundle over a curve of genus 3 with
fibre~$\CP^1$. Applying standard theory (e.g.\ \cite[p.~604--605]{GH}) we
find that each blow-up increases $b^2$ by~1 and $b^3$ by~$b^3(E)=2g(c_i)=6$
(cf.~\cite[p.~154]{Ko}), so $b^2(\oW)=4$ and $b^3(\oW)=24$. Recall that
$d(\oW)=3$, so the sublattice $L_1$ of the $K3$ lattice corresponding to
the embedding of an anticanonical $K3$ divisor in $\oW$ has rank~$r_1=1$.

As we noted earlier, the proof of Theorem~\ref{keythm} depends only on the
properties of the induced sublattices of the $K3$ lattice and since
in the present case $r_1\le 10$ the result of Theorem~\ref{emb}(a)(b)
extends as follows.
\begin{prop}\label{gen}
Let $\oW_2$ be a compact threefold of Fano type or non-symplectic type
and $L_2$ the corresponding sublattice of the $K3$ lattice with $r_2=\rk L_2$
as defined in~\S\ref{g2mfds}. If $\oW_2$ is of Fano type~$V_2$ or $\oW_2$
is of non-symplectic type $(r_2,a_2,\delta_2)$ with $r_2+a_2<20$, then a
compact irreducible $G_2$-manifold $M$ can be constructed from $\oW_2$ and
$\oW$ with
$$
b^2(M)=
\begin{cases}
3,\\
5+r_2-a_2,
\end{cases}
\qquad
b^3(M)=
\begin{cases}
b^3(V_2)-K_{V_2}^3+52,\\
96-r_2-3a_2,
\end{cases}
$$
where the top row corresponds to $\oW_2$ of Fano type and the bottom row to
$\oW_2$ of non-symplectic type.
\end{prop}

As in Theorem~\ref{emb} before, we have $n=0$ and $d_2$ even. Since in the
present case $d(\oW)$ is odd, the resulting $G_2$-manifolds have
odd~$b^2(M)$ and even $b^3$. The Betti numbers of the $G_2$-manifolds
obtainable from Proposition~\ref{gen} are:
\begin{equation*}
\begin{split}
b^2&=3, \; b^3=58,64,70,\,\ldots,\, 108,114,116,158,\\
b^2&=5, \; b^3=60,64,68,72,76,80,84,88,92,\\
b^2&=7, \; b^3=62,66,70,74,78,82,86,90,94,\\
b^2&=9, \; b^3=64,68,72,76,80,84,\\
b^2&=11,\; b^3=66,70,74,78,82,\\
b^2&=13,\; b^3=68,72,76,80,84,\\
b^2&=15,\; b^3=70,74,78,82,86,\\
b^2&=17,\; b^3=72,76\\
(b^2&,b^3)=\, (19,74)\; (21,76)\text{ and }(23,78)
\end{split}
\end{equation*}
There are 69 pairs, 10 of these overlap with the pairs $(b^2,b^3)$ given
in~\cite[figure~11.1]{Jo2}. Curiously, the overlaps occur only when $\oW_2$
is of Fano type $V_2$ with $b^3(V_2)-K_{V_2}^3$ a multiple of~4.

\subsection{Examples arising from Proposition~\ref{embed}}
\label{from.embed}
It follows from the arguments of Proposition~\ref{embed} and
Theorem~\ref{keythm} that when each of the two polarizing lattices has a
primitive embedding in $-E_8(2)\oplus H$ we can construct a compact
simply-connected $G_2$-manifold. This situation occurs when each threefold
$\oW_j$ is of Fano type $S_d\times\CP^1$, $1\le d\le 8$ or of
non-symplectic type with invariants $(10,8,0)$. The resulting $(b^2,b^3)$
are then computed as in Theorem~\ref{emb}, according to the
types of~$\oW_j$. We obtain 24 pairs
$$
(0,39+6k),\; 0\le k\le 14;\;\;
(4,43+6k_1),\; 0\le k_1\le 7\,\text{ and } (8,47).
$$

A similar possibility occurs when both $\oW_j$ are of non-symplectic type
with invariants $(r_j,a_j,\delta_j)$, in the borderline case
$r_1 + r_2 + a_1+a_2 = 20$ and $\delta_1=\delta_2=0$. The formulae of
Theorem~\ref{emb}(c) for $(b^2,b^3)$ apply in this case, giving two pairs
$(14,81), (22,89)$.

From the 26 pairs in this subsection, 4 occur in~\cite{Jo2}, 18 appeared in
the previous subsections and these 5 pairs
$$
(0,39), (0,45), (0,51), (4,49), (22,89)
$$
are new, to the authors' knowledge.

\subsection{Examples with large $\rk L_1$ and $\rk L_2 =1$}
\label{large.rank}

Finally, here are some new examples of $G_2$-manifolds with `larger' values of
$b^2$. Note that\linebreak
$L(19,1,1) \cong 2(-E_8)\oplus\mathbf 1 (-2)\oplus H$, and
$L(18,2,1)\cong 2(-E_8)\oplus\mathbf 1 (-2)\oplus\mathbf 1 (2) $, where
$\mathbf 1 (\pm 2)$ are rank 1 lattices generated by $e$ with
$e\cdot e = \pm 2$. Thus $L(19,1,1)$ and $L(18,2,1)$ can be
primitively embedded in $2(-E_8)\oplus 2H$ and the direct sum of either of
these and an even lattice $L_2$ of rank one can be primitively embedded in
the $K3$ lattice. There are 17 appropriate choices of
$\oW_2$ of Fano type arising from 17 algebraic families of Fano threefolds
$V$ with $b^2(V)=1$ and one of non-symplectic type corresponding to the
invariant lattice $L(1,1,1)$ so that Theorem~\ref{keythm}(a) applies. In
addition, we can use the threefold defined in Example~\ref{mult}. This
gives 38 instances of generalized connected sum, realizing, by calculation
similar to that in Theorem~\ref{emb}(b)(c) and Proposition~\ref{gen},
simply-connected $G_2$-manifolds $M$ with, respectively,
$$
b^2(M)=
\begin{cases}
2+r_1-a_1,\\
4+r_1-a_1,\\
5+r_1-a_1,
\end{cases}
\qquad
b^3(M)=
\begin{cases}
b^3(V_2)-K_{V_2}^3-r_1-3a_1+71,\\
111-r_1-3a_1,\\
96-r_1-3a_1.
\end{cases}
$$
Here $(r_1,a_1)$ is either $(19,1)$ or $(18,2)$.
We obtain 28 distinct pairs of Betti numbers, 13 of which are not in
the previous subsections
\begin{equation}\label{more}
\begin{split}
b^2=18,& \; b^3=69,71,83,111,153,\\
b^2=20,& \; b^3=71,73,85,95,113,155,\\
(b^2,b^3)=& \; (21,72)\text{ and }(23,74).
\end{split}
\end{equation}
For each value of $b^2$ in the above list, the respective values of $b^3$
are larger than any of those appearing in \cite[figure~11.1]{Jo2}. Thus, to
our knowledge, all of the examples of $G_2$-manifolds in~\eqref{more} are new.

\appendix
\section*{Appendix.
Invariants of Fano threefolds and non-symplectic involutions}

We give some details of classification results which were used in the
calculations of Betti numbers of $G_2$-manifolds in this paper.

Complete classification of non-singular Fano threefolds is given
in~\cite{Is,MoMu,MoMu2}; there are 105 algebraic families overall.
Here we only require the values of parameter $b^3(V)-K_V^3$ realized by Fano
threefolds~$V$. These are given in Table~\ref{fano}, sorted by the Picard
number $b^2(V)$ (as before, dots between values denote consecutive even
integers).
\begin{table}[h]
\caption{All possible values $b^2(V),b^3(V)-K_V^3$}
\begin{tabular}{|c|c|}
\hline
$b^2(V)$ & $b^3(V)-K_V^3$\\
\hline
1 & 22, 24, 26, 30, 34, 36, 40, 46, 50, 54, 64, 104\\
\hline
2 & 22, \ldots, 34, 38, 40, 42, 46, 48, 54, 56, 62\\
\hline
3 & 20, \ldots, 52\\
\hline
4 & 26, 30, \ldots, 46\\
\hline
5 & 28, 36\\
\hline
6 & 30\\
\hline
7 & 24\\
\hline
8 & 18\\
\hline
9 & 12\\
\hline
10 & 6\\
\hline
\end{tabular}
\label{fano}
\end{table}
The number of values in the second column is less than 105 as some values
$b^3(V)-K_V^3$ are realized by more than one algebraic family. If
$b^2(V) \geq 6$, then $V= S_{11-b^2(V)}\times\CP^1$ is a product with a del
Pezzo surface and $b^3(V)=-K_V^3=6(11-b^2(V))$.

We also used the classification of $K3$ surfaces with  of non-symplectic
involutions. Recall from Proposition~\ref{clas} that these are determined,
up to deformation, by the triples $(r,a,\delta)$ consisting of the rank $r$
of the invariant sublattice $L^\rho$ of Picard lattice, the number $a$ of
generators of the discriminant group $(L^\rho)^*/L^\rho$ of this
sublattice, and $\delta=\delta(L^\rho)$ the parity of involution~$\rho$.
\begin{figure}[h]
\begin{center}\begin{picture}(175,140)
\put(66,120){\circle{7}}
\put(71,118){{\tiny means $\delta=0$}}
\put(66,112){\circle*{3}}
\put(71,110){{\tiny means $\delta=1$}}

\multiput(8,10)(8,0){20}{\line(0,1){94}}
\multiput(0,18)(0,8){11}{\line(1,0){165}}
\put(0,10){\vector(0,1){100}}
\put(0,10){\vector(1,0){175}}
\put(  6,0){{\tiny $1$}}
\put( 14,0){{\tiny $2$}}
\put( 22,0){{\tiny $3$}}
\put( 30,0){{\tiny $4$}}
\put( 38,0){{\tiny $5$}}
\put( 46,0){{\tiny $6$}}
\put( 54,0){{\tiny $7$}}
\put( 62,0){{\tiny $8$}}
\put( 70,0){{\tiny $9$}}
\put( 76,0){{\tiny $10$}}
\put( 84,0){{\tiny $11$}}
\put( 92,0){{\tiny $12$}}
\put(100,0){{\tiny $13$}}
\put(108,0){{\tiny $14$}}
\put(116,0){{\tiny $15$}}
\put(124,0){{\tiny $16$}}
\put(132,0){{\tiny $17$}}
\put(140,0){{\tiny $18$}}
\put(148,0){{\tiny $19$}}
\put(156,0){{\tiny $20$}}

\put(-8, 9){{\tiny $0$}}
\put(-8,  17){{\tiny $1$}}
\put(-8, 25){{\tiny $2$}}
\put(-8, 33){{\tiny $3$}}
\put(-8, 41){{\tiny $4$}}
\put(-8, 49){{\tiny $5$}}
\put(-8, 57){{\tiny $6$}}
\put(-8, 65){{\tiny $7$}}
\put(-8, 73){{\tiny $8$}}
\put(-8, 81){{\tiny $9$}}
\put(-10, 89){{\tiny $10$}}
\put(-10, 97){{\tiny $11$}}
\put( -2,118){{\footnotesize $a$}} 
\put(186,  8){{\footnotesize $r$}} 

\put( 16,10){\circle{7}}
\put( 80,10){\circle{7}}
\put(144,10){\circle{7}}
\put( 16,26){\circle{7}}
\put( 48,26){\circle{7}}
\put( 80,26){\circle{7}}
\put(112,26){\circle{7}}
\put(144,26){\circle{7}}
\put( 48,42){\circle{7}}
\put( 80,42){\circle{7}}
\put(112,42){\circle{7}}
\put( 80,58){\circle{7}}
\put( 112,58){\circle{7}}
\put( 80,74){\circle{7}}
\put( 80,90){\circle{7}}

\put(  8,18){\circle*{3}}
\put( 24,18){\circle*{3}}
\put( 72,18){\circle*{3}}
\put( 88,18){\circle*{3}}
\put(136,18){\circle*{3}}
\put(152,18){\circle*{3}}
\put( 16,26){\circle*{3}}
\put( 32,26){\circle*{3}}
\put( 64,26){\circle*{3}}
\put( 80,26){\circle*{3}}
\put( 96,26){\circle*{3}}
\put(128,26){\circle*{3}}
\put(144,26){\circle*{3}}
\put( 24,34){\circle*{3}}
\put( 40,34){\circle*{3}}
\put( 56,34){\circle*{3}}
\put( 72,34){\circle*{3}}
\put( 88,34){\circle*{3}}
\put(104,34){\circle*{3}}
\put(120,34){\circle*{3}}
\put(136,34){\circle*{3}}
\put( 32,42){\circle*{3}}
\put( 48,42){\circle*{3}}
\put( 64,42){\circle*{3}}
\put( 80,42){\circle*{3}}
\put( 96,42){\circle*{3}}
\put(112,42){\circle*{3}}
\put(128,42){\circle*{3}}
\put( 40,50){\circle*{3}}
\put( 56,50){\circle*{3}}
\put( 72,50){\circle*{3}}
\put( 88,50){\circle*{3}}
\put(104,50){\circle*{3}}
\put(120,50){\circle*{3}}
\put( 48,58){\circle*{3}}
\put( 64,58){\circle*{3}}
\put( 80,58){\circle*{3}}
\put( 96,58){\circle*{3}}
\put(112,58){\circle*{3}}
\put( 56,66){\circle*{3}}
\put( 72,66){\circle*{3}}
\put( 88,66){\circle*{3}}
\put(104,66){\circle*{3}}
\put( 64,74){\circle*{3}}
\put( 80,74){\circle*{3}}
\put( 96,74){\circle*{3}}
\put( 72,82){\circle*{3}}
\put( 88,82){\circle*{3}}
\put( 80,90){\circle*{3}}
\put( 88,98){\circle*{3}}
\put( 96,90){\circle*{3}}
\put( 104,82){\circle*{3}}
\put( 112,74){\circle*{3}}
\put( 120,66){\circle*{3}}
\put( 128,58){\circle*{3}}
\put( 136,50){\circle*{3}}
\put( 144,42){\circle*{3}}
\put( 144,42){\circle{7}}
\put( 152,34){\circle*{3}}
\put( 160,26){\circle*{3}}
\end{picture}\end{center}

\caption{All possible invariants $(r, a, \delta)$}
\label{NiFigure}
\end{figure}
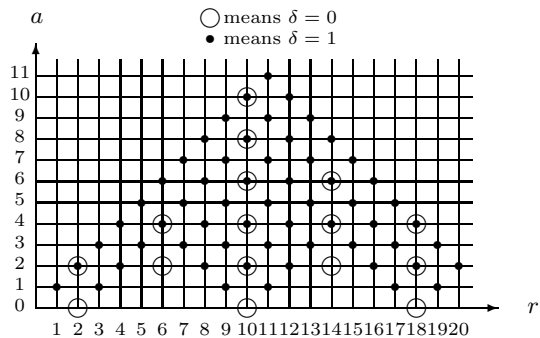
There are 75 possibilities for these invariants, these are conveniently
expressed by the Figure \ref{NiFigure}, reproduced from \cite[\S 4.5]{Ni3}.

It is also shown in~\cite[\S 4.3]{Ni3} that if $r>1$ then $L^\rho$ may be
written as $S_2\oplus(\oplus_i (-K_i))$, where $S_2$ is one of the
non-degenerate, signature $(1,1)$ lattices $H,
H(2),\mathbf{1}(2)\oplus\mathbf{1}(-2)$ and $-K_i$
are the negative-definite root lattices $-A_1=\mathbf{1}(-2), -D_{2k}, -E_7,
-E_8,-E_8(2)$. Here, as before, the lower index indicates the rank and the
notation $H(2),\, -E_8(2)$ means that the bilinear forms of~$H,\, -E_8$ are
multiplied by~2.
\vskip 2pt

{\bf Acknowledgements.}
We thank Mark Haskins for pointing out a mistake in Lemma~6.47 in~\cite{Ko};
a self-contained argument eliminating the problem now appears in the proof
of Theorem~\ref{keythm}. We also thank Mark Haskins, Viacheslav Nikulin and
Johannes Nordstr\"om for useful discussions.

\end{document}